\documentclass[11pt,thmsa]{article}
\usepackage{amsmath}
\usepackage{amssymb}

\newtheorem{Theorems1}{Theorem}[section]
\newenvironment{proof}[1][Proof]{\begin{trivlist}
\item[\hskip \labelsep {\bfseries #1}]}{\end{trivlist}}

\begin{document}

\title{Iterated Ultrapowers for the Masses}
\author{Ali Enayat \\
Department of Philosophy, Linguistics, and Theory of Science\\
University of Gothenburg, Box 200\\
SE 405 30, Gothenburg, Sweden\\
\texttt{E-mail:} \texttt{ali.enayat@gu.se} \and Matt Kaufmann \\
Department of Computer Science\\
The University of Texas at Austin\\
2317 Speedway, Stop D9500\\
Austin, TX 78712-1757 U.S.A.\\
\texttt{E-mail: kaufmann@cs.utexas.edu} \and Zachiri McKenzie \\
Department of Philosophy, Linguistics, and Theory of Science\\
University of Gothenburg, Box 200\\
SE 405 30, Gothenburg, Sweden\\
\texttt{E-mail:} \texttt{zach.mckenzie@gmail.com}}
\maketitle

\begin{abstract}
\noindent We present a novel, perspicuous framework for building iterated
ultrapowers. Furthermore, our framework naturally lends itself to the
construction of a certain type of order indiscernibles, here dubbed \textit{%
tight indiscernibles}, which are shown to provide smooth proofs of several
results in general model theory.

\begin{equation*}
\ast \ast \ast \ast \ast \ast \ast \ast \ast \ast \ast \ast \ast \ast \ast
\ast \ast \ast \ast
\end{equation*}

 \medskip

\end{abstract}

\section{Introduction and preliminaries}

\label{Intro and prelim.}

\footnotetext{%
\noindent \textbf{Acknowledgements.}~We are grateful to Jim Schmerl, John
Baldwin, and an anonymous referee for valuable feedback on earlier drafts of
this paper. Matt Kaufmann thanks the Department of Philosophy, Linguistics
and Theory of Science at the University of Gothenburg for its hospitality in
Summer 2015.}One of the central results of model theory is the celebrated
Ehrenfeucht-Mostowski theorem on the existence of models with indiscernible
sequences. All textbooks in model theory, including those of Chang \&
Keisler \cite{Chang-Keisler}, Hodges \cite{Wilfrid}, Marker \cite{Marker},
and Poizat \cite{Bruno} demonstrate this theorem using essentially the same
compactness argument that was originally presented by Ehrenfeucht and
Mostowski in their seminal 1956 paper \cite{Ehrenfeuct-Mostowski}, using
Ramsey's partition theorem. \medskip

The source of inspiration for this paper is a fundamentally different proof
of the Ehrenfeucht-Mostowski theorem that was discovered by Gaifman \cite%
{Haim-Unform} in the mid-1960s, a proof that relies on the technology of
\textit{iterated ultrapowers}, and in contrast with the usual proof, does
not invoke Ramsey's theorem. Despite the passage of several decades,
Gaifman's proof seems to be relatively unknown among logicians, perhaps due
to the forbidding technical features of the existing constructions of
iterated ultrapowers in the literature.\footnote{%
In contrast, iterated ultrapowers are well-known to specialists in the
`large cardinals' area of set theory, albeit in a very specialized form,
where (a) the ultrafilters at work have extra combinatorial features, and
(b) the iteration is always carried out along a well-ordering (as opposed to
a linear ordering).} Here we attempt to remedy this situation by presenting
a streamlined account that is sufficiently elementary to be accessible to
logicians familiar with the rudiments of model theory. \medskip

Our exposition also incorporates novel technical features: we bypass the
usual method of building iterated ultrapowers as direct limits, and instead
build them in the guise of \emph{dimensional }(\emph{Skolem})\emph{\
ultrapowers}, in a natural manner reminiscent of the usual construction of
(Skolem) ultrapowers. We also isolate the key notion of \emph{tight
indiscernibles} to describe a special and useful kind of order
indiscernibles that naturally arise from dimensional ultrapowers. In the
interest of balancing clarity with succinctness, we have opted for a
tutorial style: plenty of motivation is offered, but some of the proofs are
left in the form of exercises for the reader, and the solutions are
collected in the appendix. \medskip

In the rest of this section we review some preliminaries, introduce the key
notion of tight indiscernibles\emph{,} and state a corresponding existence
theorem (Theorem 1.3). Then in Section 2 we utilize tight indiscernibles to
prove some results in general model theory, including a refinement of the
Ehrenfeucht-Mostowski Theorem (Theorem 2.4). Having thus demonstrated the
utility of tight indiscernibles in Section 2, we show how to construct them
with the help of dimensional ultrapowers in Section 3. Section 4 presents a
brief historical background of iterated ultrapowers, with pointers to the
rich literature of the subject.\medskip

\noindent \textbf{1.1.}~\textbf{Definitions, notations and conventions.}~All
the structures considered in this paper are first order structures; we
follow the convention of using $M$, $M^{\ast },$ $M_{0}$, etc.~to denote
(respectively) the universes of discourse of structures $\mathcal{M}$, $%
\mathcal{M}^{\ast },$ $\mathcal{M}_{0},$ etc. We assume the Axiom of Choice
in the metatheory, and use $\omega $ for the set of non-negative integers.
\smallskip

\noindent \textbf{(a)} Given a structure $\mathcal{M}$, $\mathcal{L}(%
\mathcal{M})$ is the \emph{language} of $\mathcal{M}$, and the \textit{%
cardinality} of $\mathcal{M}$ refers to the cardinality of $M$.\smallskip

\noindent \textbf{(b)} $\mathcal{M}_{1}$ is an \emph{expansion} of $\mathcal{%
M}$ if $M_{1}=M$, $\mathcal{L}(\mathcal{M}_{1})\supseteq \mathcal{L}(%
\mathcal{M})$, and $\mathcal{M}_{1}$ and $\mathcal{M}$ give the same
interpretation to every symbol in $\mathcal{L}(\mathcal{M})$. $\mathcal{M}%
^{+}$ is the expansion $\left( \mathcal{M},m\right) _{m\in M}$ of $\mathcal{M%
}$ (we follow the practice of confusing the constant symbol in $\mathcal{L}(%
\mathcal{M}^{+})$ interpreted by $m$ with $m$ in the interest of a lighter
notation). $\mathcal{M}^{\#}$ is a \emph{full expansion} of $\mathcal{M}$%
\emph{\ }if $\mathcal{M}^{\#}=(\mathcal{M},R)_{R\in \mathcal{F}}$, where $%
\mathcal{F}=\left\{ R\subseteq M^{n}:n\in \omega \right\} .$\smallskip

\noindent \textbf{(c)} For a formula $\varphi (x_{1},\ldots ,x_{k})$ with
free variables implicitly ordered as shown, we write $\varphi ^{\mathcal{M}}$
for $\left\{ \left\langle m_{1},\ldots ,m_{k}\right\rangle \in M^{k}:%
\mathcal{M\models \varphi }\left( m_{1},\ldots ,m_{k}\right) \right\} $.
\smallskip

\noindent \textbf{(d)} Suppose $X\subseteq M^{k},$ where $k$ is a positive
integer. $X$ is $\mathcal{M}$-\emph{definable} if $X=\varphi ^{\mathcal{M}}$
for some $\mathcal{L}(\mathcal{M})$-formula. $X$ is \emph{parametrically} $%
\mathcal{M}$-\emph{definable} if $X=\varphi ^{\mathcal{M}^{+}}$ for some $%
\mathcal{L}(\mathcal{M}^{+})$-formula.\smallskip

\noindent \textbf{(e) }A \emph{parametrically} $\mathcal{M}$-\emph{definable}
function is a function $f:M^{k}\rightarrow M$ (where $k$ is a positive
integer$)$ such that the graph of $f$ is parametrically $\mathcal{M}$%
-definable. If $\mathcal{M}^{\ast }$ is an elementary extension of $\mathcal{%
M}$ (written \textquotedblleft $\mathcal{M}^{\ast }\succ \mathcal{M}$" or
\textquotedblleft $\mathcal{M}\prec \mathcal{M}^{\ast }$"), then any such $f$
extends naturally to a parametrically $\mathcal{M}^{\ast }$-definable
function according to the same definition; we may also denote this extension
as $f$.\smallskip

\noindent \textbf{(f) }$\mathcal{M}$ has \emph{definable Skolem functions}
if for every $\mathcal{L}(\mathcal{M})$-formula $\varphi (x,y_{1},\ldots
,y_{k})$, whose free variable(s) include a distinguished free variable $x$
and whose other free variables (if any) are $y_{1},\ldots ,y_{k}$, there is
an $\mathcal{M}$-definable function $f$ such that (abusing notation
slightly):

\begin{center}
$\mathcal{M}\models \forall y_{1}\ldots \forall y_{k}\left( \exists x\
\varphi (x,y_{1},\ldots ,y_{k})\rightarrow \varphi (f(y_{1},\ldots
,y_{k}),y_{1},\ldots ,y_{k})\right) .$
\end{center}

\noindent \textbf{(g)} Given a structure $\mathcal{M}$ with definable Skolem
functions and a subset $S$ of $M$, there is a least elementary substructure $%
\mathcal{M}_{S}$ of $\mathcal{M}$ that contains $S$, whose universe is the
set of all applications of $\mathcal{M}$-definable functions to tuples from $%
S$. $\mathcal{M}_{S}$ is called the \emph{submodel of }$\mathcal{M}$\emph{\
generated by} $S$. \smallskip

\noindent \textbf{(h)} Suppose $\mathcal{M}^{\ast }\succ \mathcal{M}$ and $%
(I,<)$ is a linear order with $I\subseteq M^{\ast }\backslash M$ (we often
omit explicit mention of the order, $<)$. We say that $(I,<)$ \emph{forms a
set of order indiscernibles over} $M$ if for any $\mathcal{L}(\mathcal{M})$%
-formula $\varphi (x_{1},...,x_{n},y_{1},\ldots ,y_{k})$, any sequence $%
m_{1},...,m_{n}$ from $M$, and any two sequences $i_{1}<\ldots <i_{k}$ and $%
j_{1}<\ldots <j_{k}$ from $I$, we have:

\begin{center}
$\left( \mathcal{M}^{\ast },m\right) _{m\in M}\models \varphi
(m_{1},...,m_{n},i_{1},\ldots ,i_{k})\leftrightarrow \varphi
(m_{1},...,m_{n},j_{1},\ldots ,j_{k}).$
\end{center}

With the above preliminaries in place, we are ready for a definition that
plays a central role in this paper.\medskip

\noindent \textbf{1.2.}~\textbf{Definition.} Suppose $\mathcal{M}$ has
definable Skolem functions, $\mathcal{M}^{\ast }\succ \mathcal{M}$ and $%
(I,<) $ is a linear order with $I\subseteq M^{\ast }\backslash M.$ $(I,<)$
\emph{forms a set of tight indiscernibles generating} $\mathcal{M}^{\ast }$
\emph{over} $\mathcal{M}$ if the following three properties hold. (Note: For
brevity we will say that $I$ \emph{forms a set of tight indiscernibles over}
$\mathcal{M}$, when the order and $\mathcal{M}^{\ast }$ are clear from the
context.)\smallskip

\noindent \textbf{(1)} $(I,<)$ is a set of order indiscernibles over $%
\mathcal{M}$.\smallskip

\noindent \textbf{(2)} $\mathcal{M}^{\ast }$ is generated by $M\cup I$,
i.e., every element of $\mathcal{M}^{\ast }$ is of the form $f(i_{1},\ldots
,i_{k})$ for some $i_{1},\ldots ,i_{k}$ from $I$ and some parametrically $%
\mathcal{M}$-definable function $f$.\smallskip

\noindent \textbf{(3)} For all $i_{1}<\ldots <i_{k}<j_{1}<\ldots <j_{k}$
from $I$ and every parametrically $\mathcal{M}$-definable function $f$, if $%
f(i_{1},\ldots ,i_{k})=f(j_{1},\ldots ,j_{k})$ then this common value is in $%
M$.\medskip

\noindent \textbf{1.2.1.}~\textbf{Remark.~}Property (3) of tight
indiscernibles is not implied by the other two properties: Start with any
structure $\mathcal{M}$ with definable Skolem functions that also has a
definable bijection $\pi $ from $M^{2}$ to $M$. Use the
Ehrenfeucht-Mostowski theorem to build a proper chain $\mathcal{M}\prec
\mathcal{M}_{1}\prec \mathcal{M}_{2}$ such that for some linear order $%
\left( I_{0},<\right) ,$ $I_{0}\subseteq M_{2}\backslash M_{1},$ and $\left(
I_{0},<\right) $ is a set of order indiscernibles in $\mathcal{M}_{2}$ over $%
\mathcal{M}_{1}.$ Then fix some $a\in M_{1}\backslash M$ and let $I=\{\pi
(a,i):i\in I_{0}\},$ with the obvious order inherited from $I$. Finally, let
$\mathcal{M}^{\ast }$ be the elementary submodel of $\mathcal{M}_{2}$
generated by $M\cup I.$ In this example properties (1) and (2) of Definition
1.2 are satisfied, but not property (3), since $a\notin M$ and yet $f(i)=a$
for every $i\in I$, where $f$ is the definable function given by $f(x)=y$
iff $y$ is the unique element $y$ such that $\exists z\ \pi (y,z)=x.$\medskip

\noindent \textbf{1.3.}~\textbf{Theorem~}(Existence of Tight Indiscernibles)
\textit{Every infinite structure} $\mathcal{M}$ \textit{has an expansion} $%
\mathcal{M}_{1}$ \textit{with definable Skolem functions} \textit{such that
for} \textit{any} \textit{ordered set} $(I,<)$, \textit{with} $I$ \textit{%
disjoint from} $M_{1}$, \textit{there is} $\mathcal{M}^{\ast }\succ \mathcal{%
M}_{1}$ \textit{such that} $(I,<)$ \textit{forms a set of tight
indiscernibles generating} $\mathcal{M}^{\ast }$ \textit{over} $\mathcal{M}%
_{1}$. \textit{Moreover,} $\left\vert \mathcal{L}(\mathcal{M}%
_{1})\right\vert =\max \left\{ \left\vert \mathcal{L}(\mathcal{M}%
)\right\vert ,\aleph _{0}\right\} $.\medskip

\noindent The proof of Theorem 1.3 will be presented in Section 3, using the
technology of dimensional ultrapowers for all infinite structures $\mathcal{M%
}$. However, if both\textit{\ }$\mathcal{L}(\mathcal{M})$ and $M$ are
countable, then a straightforward proof can be given by fine-tuning the
Ehrenfeucht-Mostowski method of building order indiscernibles:\medskip

\noindent \textbf{Proof of Theorem 1.3 when} $\left\vert \mathcal{L}(%
\mathcal{M})\right\vert \leq \left\vert M\right\vert =\aleph _{0}.$ Let $%
\mathcal{M}_{1}=(\mathcal{M},\vartriangleleft )$, where $\vartriangleleft $
is an ordering of $M$ of order type $\omega $. Note that $\mathcal{M}_{1}$
has definable Skolem functions. Let $\left\langle \varphi _{i}:i\in \omega
\right\rangle $ be an enumeration $\mathcal{L}(\mathcal{M}_{1}^{+})$%
-formulae with at least one free variable, and suppose the set of free
variables of $\varphi _{i}$ is $\left\{ x_{1},...,x_{n_{i}}\right\} .$ Use
Ramsey's partition theorem to construct a decreasing sequence $\left\langle
X_{i}:i\in \omega \right\rangle $ of infinite subsets of $M$ such that for
all $i\in \omega $ either $\varphi _{i}$ or $\lnot \varphi _{i}$ holds for
all increasing $a_{1}\vartriangleleft \ldots \vartriangleleft a_{n_{i}}$
from $X_{i}$. Fix an ordered set $(I,<)$ disjoint from $M.$ Let $T$ be the
corresponding Ehrenfeucht-Mostwoski blueprint, i.e., the result of
augmenting the elementary diagram of $\mathcal{M}_{1}$ with all sentences $%
\varphi (i_{1},\ldots ,i_{n})$, where $i_{1}<\ldots <i_{n}$ are elements of $%
I$, $\varphi $ is allowed to have parameters from $M,$ and for some $k\in
\omega $ and all $a_{1}\vartriangleleft \ldots \vartriangleleft a_{n}$ from $%
X_{k}$, $\varphi (a_{1},\ldots ,a_{n})$ holds in $\mathcal{M}$. $T$ is
easily seen to be consistent by compactness considerations. We obtain the
desired $\mathcal{M}^{\ast }\succ \mathcal{M}_{1}$ in which $I$ forms a set
of tight indiscernibles by starting with any model of $T$ and taking the
submodel generated by $M$ and $I$. Now suppose that in $\mathcal{M}^{\ast }$
we have strictly increasing sequences from $I$, $c_{1}<\ldots <c_{k}$ and $%
d_{1}<\ldots <d_{k}$, with $c_{k}<d_{1}$ such that:\smallskip

\noindent $\lbrack \ast \rbrack$\quad $f(c_{1},\ldots ,c_{k})=f(d_{1},\ldots
,d_{k})$.\smallskip

\noindent We shall show that this (common) value is in $M$. For some $X_{i}$%
, we have $f(a_{1},\ldots ,a_{k})=f(b_{1},\ldots ,b_{k})$ for all $%
a_{1}<\ldots <a_{k}<b_{1}<\ldots <b_{k}$ from $X_{i}$ (since if instead the
equality were false on $X_{i}$, $[\ast ]$ would fail). Let $m\in M$ be the
common value of all such $f\left( a_{1},\ldots ,a_{k}\right) $. Then for
some $X_{j}$, $m=f\left( a_{1},\ldots ,a_{k}\right) $ for all $%
a_{1}\vartriangleleft \ldots \vartriangleleft a_{k}$ from $X_{j}$, since
otherwise this equality is false for all $a_{1}\vartriangleleft \ldots
\vartriangleleft a_{k}$ from some $X_{j}$, which is impossible by
considering what holds in $X_{s}$ for $s=\max \{i,j\}$. So $f(c_{1},\ldots
,c_{k})=m\in M$. \hfill $\square $\bigskip

\section{Tight indiscernibles at work}

\label{At work}

In this section we show that the existence of tight indiscernibles provides
straightforward proofs of several results of general model theory. The
following lemma earns us a key strengthening of property (3) of tight
indiscernibles.\footnote{%
An abstract formalization of this lemma has been verified with the ACL2
theorem prover, cf.~\cite{Matt's ACL2 Proof}.}\medskip

\noindent \textbf{2.1.~Tightness Lemma.}~\textit{Suppose that} $(I,<)$
\textit{is an infinite linear order and }$I$\textit{\ forms a set of tight
indiscernibles generating} $\mathcal{M}^{\ast }$ \textit{over }$\mathcal{M}$%
\textit{. Also suppose that }$f$ \textit{and} $g$ \textit{are parametrically
}$\mathcal{M}$\textit{-definable functions and there are disjoint subsets} $%
S_{0}=\left\{ i_{1},\ldots ,i_{p}\right\} $ \textit{and} $S_{1}=\left\{
j_{1},\ldots ,j_{q}\right\} $ \textit{of} $I$ \textit{such that }$%
f\left(i_{1},\ldots ,i_{p}\right) =g\left( j_{1},\ldots ,j_{q}\right) $%
\textit{. Then} $f\left( i_{1},\ldots ,i_{p}\right) \in M$.\medskip

\noindent \textbf{Proof.}~We may assume (by tweaking $f$ and $g$) that $%
i_{1}<\ldots <i_{p}$ and $j_{1}<\ldots <j_{q}.$ Let the (disjoint) union $%
S_{0}\cup S_{1}$ be ordered as $k_{1}<\ldots <k_{r}$, where $r=p+q$. Since $%
(I,<)$ is infinite, it contains a strictly increasing or strictly decreasing
$\omega $-sequence (by Ramsey's theorem, or as one of us was told years ago:
the order is either well-founded or it's not). Without loss of generality
assume there is a decreasing such sequence, as the other case is entirely
analogous. By indiscernibility, we may assume that the elements $%
k_{1},\ldots ,k_{r}$ are contained in this sequence; hence we may choose $%
k_{1}^{\prime }<\ldots <k_{r}^{\prime }$ such that $k_{r}^{\prime }<k_{1}.$%
\medskip

Our plan is to shift $k_{1}$ way to the left, then shift $k_{2}$ way to the
left (but still to the right of the new $k_{1}$), and so on, so that the new
$k_{r}$ is shifted to the left of the original $k_1$. So for each natural
number $j\leq r$, define the sequence $\left\langle n_{\{j,i\}}:1\leq
i<r\right\rangle $ as follows: $n_{\{j,i\}}=k_{i}^{\prime }$ for $i\leq j,$
and $n_{\{j,i\}}=k_{i}$ for $i>j$. In particular, we have

\begin{center}
$\left\langle n_{\{0,i\}}:1\leq i\leq r\right\rangle =\left\langle
k_{i}:1\leq i\leq r\right\rangle $ and

$\left\langle n_{\{r,i\}}:1\leq i\leq r\right\rangle =\left\langle
k_{i}^{^{\prime }}:1\leq i\leq r\right\rangle $.
\end{center}

\noindent Let $I_{0}$, $J_{0}$ be the indices $i$ for which $k_{i}\in S_{0}$
or $k_{i}\in S_{1}$, respectively. Let $m=f(i_{1},\ldots ,i_{p})$. Then
indiscernibility and a straightforward induction on $j$ (and a slight but
clear abuse of notation) establish the following. (\textit{Hint}:
disjointness of $S_{0}$ and $S_{1}$ guarantee that when moving from $j$ to $%
j+1$ for the induction step, either the application of $f$ or the
application of $g$ remains unchanged.)

\begin{center}
$m=f\left( n_{\{j,s\}}:s\in I_{0}\right) =g\left( n_{\{j,s\}}:s\in
J_{0}\right) .$
\end{center}

\noindent In particular, instantiating with $j=0$ and $j=r$ we see that

\begin{center}
$m=f(n_{\{0,s\}}:s\in I_{0})=f(n_{\{r,s\}}:s\in I_{0})$.
\end{center}

\noindent But $n_{\{r,s_{1}\}}<n_{\{0,s_{2}\}}$ for all $s_{1}$ and $s_{2}$;
so property (3) of tight indiscernibles immediately implies that $m\in M$.
\hfill $\square $\medskip

\noindent \textbf{2.1.1.}~\textbf{Remark.}~After reading our paper, Jim
Schmerl observed that Lemma 2.1 holds for all ordered sets $I$\ whose size
exceeds $p+q$.~ We leave the proof to the reader. \textit{Hint}: Assuming
without loss of generality that $p\leq q$, first observe that by shifting
the given arguments of $f$ and $g$, we see that the value of $f$ is
preserved when the only change is to shift just one of the arguments without
changing their ordering.~ Now successively shift all arguments of $f$ to the
left with maximum argument $u$, and then again to the right so that all
arguments are beyond $u$.~ The definition of tightness now applies.\medskip

We are now ready to use tight indiscernibles to establish the following
theorems. \medskip

\noindent \textbf{2.2.}~\textbf{Theorems.}~(Applications of Tight
Indiscernibles) \textit{Let} $\mathcal{M}$ \textit{be any infinite structure}%
.\textit{\smallskip }

\noindent \textbf{(A)} \textit{There is a proper elementary extension} $%
\mathcal{M}^{\ast}$ \textit{of} $\mathcal{M}$ \textit{such that for some
automorphism} $\alpha $ \textit{of }$\mathcal{M}^{\ast }$, $\alpha $ \textit{%
and all its finite iterates} $\alpha ^{n}$ \textit{fix each element of }$M$
\textit{and move every element in} $M^{\ast }\backslash M.$ \textit{Indeed} $%
\mathcal{M}^{\ast}$ \textit{can be arranged to be of any cardinality} $%
\kappa \geq \max \{\left\vert M\right\vert ,\left\vert \mathcal{L}(\mathcal{M%
})\right\vert \}.$\medskip

\noindent \textbf{(B)} \textit{There is a family of proper elementary
extensions} $\{\mathcal{M}_{i}:i\in \omega \}$ \textit{of} $\mathcal{M}$%
\textit{\ whose intersection is }$M$, \textit{such that for all} $i$, $%
\mathcal{M}_{i}$ \textit{is a} \textit{proper elementary extension of }$%
\mathcal{M}_{i+1}$.\medskip

\noindent \textbf{(C)} \textit{Let} $I$ \textit{be any set. Then there is a
family} $\{\mathcal{M}_{S}:S\subseteq I\}$ \textit{of proper elementary
extensions of} $\mathcal{M}$ \textit{such that} $\mathcal{M}_{S_{0}}\prec
\mathcal{M}_{S_{1}}$\textit{\ whenever} $S_{0}\subseteq S_{1}\subseteq I$
\textit{with the additional property that }$M_{S_{0}}\cap M_{S_{1}}=M$
\textit{for all disjoint} $S_{0},$ $S_{1}$.\footnote{%
The proof of (C) shows that we may conclude a bit more: for each $S\subseteq
I$, $S$ forms a set of tight indiscernibles in $\mathcal{M}_{S}$. However,
we find the theorem to be of interest even without that embellishment.}%
\medskip

\noindent \textbf{Proof.} Observe that by taking reducts back to $\mathcal{L}%
(\mathcal{M})$ it suffices to demonstrate Theorem 2.2 for some \emph{%
expansion} of $\mathcal{M}$. To this end, we may assume without loss of
generality that $\mathcal{M}$ is replaced by the expansion $\mathcal{M}_{1}$
of Theorem 1.3. We begin each proof by applying Theorem 1.3 to obtain an
elementary extension $\mathcal{M}^{\ast }$ of $\mathcal{M}$ with a set of
tight indiscernibles generating $\mathcal{M}^{\ast }$ over $\mathcal{M}$. We
now apply the Tightness Lemma to obtain each result in turn, as
follows.\medskip

\noindent \textbf{Proof of (A).}~We first prove the theorem when $\kappa
=\max \{\left\vert M\right\vert ,\left\vert \mathcal{L}(\mathcal{M}%
)\right\vert \}$ and then we will explain how to handle larger values of $%
\kappa $ by a minor variation. Let $(I,<)$ be the ordered set $\mathbb{Z}$\
of integers. Without loss of generality assume that $\mathcal{M}$ and $I$
are disjoint (else proceed below using an isomorphic copy of $I$). We define
an automorphism $\alpha $ of $\mathcal{M}^{\ast }$ as follows: given $x\in
M^{\ast }$, by property (2) of tight indiscernibles we may write $%
x=f(i_{1},\ldots ,i_{n})$ for some $i_{1}<\ldots <i_{n}$ from $I$ and some
parametrically $\mathcal{M}$-definable function $f$. Then we define $\alpha
(x)=f(i_{1}+1,\ldots ,i_{n}+1)$. Then $\alpha $ is well-defined and is an
automorphism of $\mathcal{M}^{\ast }$, by properties (1) and (2) of tight
indiscernibles. Moreover, $\alpha $ fixes each $m\in M$ since $%
m=f_{m}(i_{1}) $, where $i_{1}\in I$ and $f_{m}$ is the constant function
with range $\{m\}$ and therefore $\alpha (m)=f_{m}(i_{1}+1)=m.$ Finally
suppose $x=\alpha (x)$ for $x$ as above; we show $x\in M$. Let $\alpha ^{k}$
be the $k$-fold composition of $\alpha $. A trivial induction shows that for
all positive integers $k$, $\alpha ^{k}$ is an automorphism of $\mathcal{M}%
^{\ast }$ and $\alpha ^{k}(x)=f(i_{1}+k,\ldots ,i_{n}+k)$. Taking $%
k=i_{n}-i_{1}+1$, we have $f(i_{1},\ldots ,i_{n})=f(i_{1}+k,\ldots ,i_{n}+k)$%
, where $i_{n}<i_{1}+k$; so by property (3) of tight indiscernibles, $x\in M$%
. It should be clear that every finite iterate of $\alpha $ pointwise fixes $%
M$ and moves every element of $M^{\ast }\backslash M.$ Finally, to handle
the case when $\kappa >\{\left\vert M\right\vert ,\left\vert \mathcal{L}(%
\mathcal{M})\right\vert \}$, let $(I,<)$ be the lexicographically ordered
set $\kappa \times \mathbb{Z}$, where $\kappa $ carries its natural order,
and take advantage of the automorphism $(\gamma ,i)\mapsto (\gamma ,i+1)$ of
$\kappa \times \mathbb{Z}$. \hfill $\square $\medskip

\noindent \textbf{Proof of (B).}~Let $(I,<)$ be the ordered set of
non-negative integers, and for each $i$ let $\mathcal{M}_{i}$ be the
submodel of $\mathcal{M}^{\ast }$ generated by $M\cup \{n\in I:n\geq i\}$.
(Thus, $\mathcal{M}_{0}=\mathcal{M}^{\ast }$.) Property (3) clearly implies
that $i\notin M_{i+1}$; hence $\mathcal{M}_{i+1}$ is a proper subset of $%
\mathcal{M}_{i}$. It remains to show that every element of the intersection
of the $\mathcal{M}_{i}$ is an element of $M$. So consider an arbitrary
element $m=f(i_{1},\ldots ,i_{k})$ of that intersection, where $i_{1}<\ldots
<i_{k}$ from $I$ and $f$ is an $\mathcal{M}$-definable function. Let $%
p=i_{k}+1$. Since $m\in M_{p}$, we may choose $j_{0}<\ldots <j_{n}$ with $%
p\leq j_{0}$ such that $m=g(j_{0},\ldots ,j_{n})$, for some $\mathcal{M}$%
-definable function $g$. Since $i_{k}<j_{0}$, we have $m\in M$ by the
Tightness Lemma. \hfill $\square $\medskip

\noindent \textbf{Proof of (C).}~Without loss of generality, assume that $I$
is ordered without endpoints, as we can always restrict to the model
generated by the original set $I$ (which can thus even be supplied with a
given order). For each $S\subseteq I$, let $\mathcal{M}_{S}$ be the submodel
of $\mathcal{M}^{\ast }$ generated by $S$. So $\mathcal{M}_{S_{0}}\prec
\mathcal{M}_{S_{1}}$ whenever $S_{0}\subseteq S_{1}\subseteq I$. Now suppose
that $S_{0}$ and $S_{1}$ are disjoint subsets of $I$, and suppose $m\in
M_{S_{0}}\cap M_{S_{1}}$; it remains only to prove that $m\in M$. But since $%
S_{0}$ and $S_{1}$ are disjoint and generate $\mathcal{M}_{S_{0}}$ and $%
\mathcal{M}_{S_{1}}$ respectively, this is immediate by the Tightness Lemma.
\hfill $\square $\medskip

\noindent \textbf{2.2.1.}~\textbf{Remark.}~As pointed out to us by Jim
Schmerl, Theorem 2.2(A) can be derived from a powerful result due to Duby
\cite{Duby} who proved that every structure $\mathcal{M}$ has an elementary
extension that carries an automorphism $\alpha $ such that $\alpha $ and all
its finite iterates $\alpha ^{n}$ are `maximal automorphisms', i.e., they
move every nonalgebraic element of $\mathcal{M}$ (an element $m_{0}$ of
structure $\mathcal{M}$ is algebraic in $\mathcal{M}$ iff there is some
unary $\mathcal{L}(\mathcal{M})$-formula $\varphi (x)$ such that $\varphi ^{%
\mathcal{M}}$ is finite and contains $m_{0}$).\footnote{%
Duby's result is a generalization of a key result of K\"{o}rner \cite%
{Koerner}, who proved that if both $\mathcal{L}(\mathcal{M)}$ and $M$ are
countable, then $\mathcal{M}$ has a countable elementary extension that
carries a maximal automorphism. K\"{o}rner's result, in turn, was inspired
by and generalizes a theorem due to the joint work of Kaye, Kossak, and
Kotlarski \cite{Richard-Roman-Henryk} that states that every countable model
of \textrm{PA} has a countable elementary extension that carries a maximal
automorphism (indeed it is shown in \cite{Richard-Roman-Henryk} that
countable recursively saturated models of $\mathrm{PA}$ with maximal
automorphisms are precisely the countable arithmetically saturated ones).}
Also, Theorem 2.2(B) appears as Exercise 3.3.9 of Chang \& Keisler's
textbook \cite{Chang-Keisler}, albeit a rather challenging one.\footnote{%
There are two reasons we included Theorem 2.2(B): (1) We know of at least
two competent logicians (a set theorist and a model theorist) who were
stumped by this exercise after assigning it to their students in a graduate
model theory course since they assumed when assigning the exercise that the
result follows easily from the usual formulation of the
Ehrenfeucht-Mostowski theorem on the existence of order indiscernibles; (2)
John Baldwin has informed us that a special case of Theorem 2.2(B) is
contained in a result of Shelah that appears as [Ba, Lemma 7.5].}\medskip

\noindent \textbf{2.3.}~\textbf{Corollary.}~\textit{If }$T$\textit{\ is a
theory that has an infinite algebraic model }(\textit{i.e., an infinite
model in which every element is algebraic})\textit{, then }$T$\textit{\ has
a nonalgebraic model of every cardinality} $\kappa \geq \max \{\aleph
_{0},\left\vert \mathcal{L}(\mathcal{M})\right\vert \}$ \textit{that carries
an automorphism }$\alpha $ \textit{such that }$\alpha $ \textit{and all its
finite iterates} $\alpha ^{n}$ \textit{are maximal automorphisms. }\medskip

\noindent \textbf{Proof.}~Let $\mathcal{M}$ be an infinite algebraic model
of $T$, and then use Theorem 2.2(A). \hfill $\square $\medskip

\noindent \textbf{2.3.1.}~\textbf{Remark.}~Every consistent theory $T$\ with
definable Skolem functions has an algebraic model since if $\mathcal{M}_{1}$
is a model of $T$ with definable Skolem functions, then by Tarski's test for
elementarity the submodel $\mathcal{M}$ of $\mathcal{M}_{1}$ whose universe
consists of the pointwise definable elements is an elementary submodel of $%
\mathcal{M}_{1}$, and thus is an algebraic model of $T.$ Corollary 2.3 can
also be derived from Duby's theorem mentioned in Remark 2.2.1.\medskip

Let $\mathrm{Aut}(\mathcal{S})$\textit{\ }be the automorphism group of the
structure $\mathcal{S}$. The usual proof of the Ehrenfeucht-Mostowski
theorem on the existence of order indiscernibles makes it clear that given
an infinite structure $\mathcal{M}$\ and a linear order $(I,<)$ there is an
elementary extension $\mathcal{M}^{\ast }$ of $\mathcal{M}$\ and a group
embedding $\alpha \mapsto \widehat{\alpha }$ from $\mathrm{Aut}(I,<)$\ into $%
\mathrm{Aut}(\mathcal{M}^{\ast })$ such that $\widehat{\alpha }$\ fixes
every element of $M$. The following theorem refines this embedding
result.\medskip

\noindent \textbf{2.4.}~\textbf{Theorem.}~\textit{Suppose }$\mathcal{M}$%
\textit{\ is an infinite structure. Given any linear order }$(I,<)$\textit{\
there is a proper elementary extension }$\mathcal{M}^{\ast }$ \textit{of }$%
\mathcal{M}$\textit{\ and a group embedding }$\alpha \mapsto \widehat{\alpha}
$ \textit{from }$\mathrm{Aut}(I,<)$\textit{\ into }$\mathrm{Aut}(\mathcal{M}%
^{\ast })$ \textit{such that} $\widehat{\alpha }$\textit{\ \textbf{moves
every element of} }$M^{\ast }\backslash M$ \textit{and fixes every element of%
} $M$ \textit{whenever }$\alpha $ \textit{is fixed point free}. \medskip

\noindent \textbf{Proof.}~The proof is an elaboration of the proof of
Theorem 2.2(A). Use Theorem 1.3 to obtain $\mathcal{M}^{\ast }\succ \mathcal{%
M}$ such that $(I,<)$ is a set of tight indiscernibles generating $\mathcal{M%
}^{\ast }$ over $\mathcal{M}$. Given $\alpha \in \mathrm{Aut}(I,<)$ let $%
\widehat{\alpha }\in \mathrm{Aut}(\mathcal{M}^{\ast })$ be given by:

\begin{center}
$\widehat{\alpha }(f(i_{1},\ldots ,i_{k}))=f\left( \alpha (i_{1}),\ldots
,\alpha (i_{k})\right) .$
\end{center}

\noindent The key observation is that for any fixed point free automorphism $%
\alpha$ of $(I,<)$ and any finite $i_{1},\ldots ,i_{k}\in I$, there is some $%
m\in \omega $ for which $\alpha^{m}$ has the property:

\begin{center}
$\left\{ \alpha^{m}(i_{1}),\ldots ,\alpha ^{m}(i_{k})\right\} \cap
\{i_{1},\ldots ,i_{k}\}=\varnothing $.
\end{center}

\noindent To see this, it suffices to note that if $\alpha $ is fixed point
free, then given $i$ and $j$ in $I$, if $\alpha ^{m_{0}}(i)=j$ for some $%
m_{0}$, then $\alpha ^{m}(i)\neq j$ for all $m>m_{0}$ (because for any $i\in
I$ either $\alpha ^{m}(i)<\alpha ^{n}(i)$ whenever $m<n\in \omega ;$ or $%
\alpha ^{m}(i)>\alpha ^{n}(i)$ whenever $m<n\in \omega $).

If $\widehat{\alpha }(f(i_{1},\ldots ,i_{k}))=f(i_{1},\ldots ,i_{k})$ for
some $f$ and some $i_{1},\ldots ,i_{n}$, then by indiscernibility $f\left(
\alpha ^{n}(i_{1}),\ldots ,\alpha ^{n}(i_{k})\right) =f(i_{1},\ldots ,i_{k})$
for all $n\in \omega $. The above observation now allows us to apply the
Tightness Lemma to conclude $f(i_{1},\ldots ,i_{k})\in M.$\footnote{%
The key observation used in the proof is a special case of P.~M.~Neumann's
Separation Lemma \cite{Kaye and Macpherson}, which states that if $G$ is a
group of permutations of an infinite set $A$ with the property that for each
$a\in A$ the $G$-orbit $\left\{ \alpha (a):\alpha \in G\right\} $ of each $%
a\in A$ is infinite, then for all finite subsets $X$ and $Y$ of $A$, there
is some $\alpha \in G$ such that $\alpha (X)\cap Y=\varnothing .$ Note that
if $\alpha $ is a fixed point free automorphism of a linear order $I$, and $%
G=\{\alpha ^{m}:m\in \mathbb{Z}\}$, then the $G$-orbit of each $i\in I$ is
infinite.}\hfill $\square $

\bigskip

\section{Dimensional Skolem ultrapowers}

\label{Constructing}

In this section we construct \emph{dimensional }(\emph{Skolem})\emph{\
ultrapowers}, typically called \emph{iterated ultrapowers} in the
literature, and show how they give rise to tight indiscernibles. Iterated
ultrapowers are typically built by means of a direct limit construction;
however, we take a different route here that is more in tune with the
construction of ordinary ultrapowers. A high-level description of
dimensional ultrapowers will be given in Subsection \ref{Motivation and
overview} to give a context for the nitty-gritty details of the remaining
sections. Finite dimensional ultrafilters are treated in Subsection \ref%
{Finite dimensional ultrafilters}. Then in Subsection \ref{Warmup on two
dim. ultrapowers} we focus on two-dimensional ultrapowers; and in Subsection %
\ref{Dimensional Ultrapowers} we build $I$-dimensional ultrapowers for any
linear order $I$ and give a proof of Theorem 1.3.

\subsection{Motivation and overview}

\label{Motivation and overview}

Fix a structure $\mathcal{M}$ with definable Skolem functions; we will
consider non-principal ultrafilters $\mathcal{U}$ on the Boolean algebra of
parametrically $\mathcal{M}$-definable subsets of $M$. The \emph{Skolem
ultrapower} of $\mathcal{M}$ modulo $\mathcal{U}$, here denoted $\mathrm{Ult}%
(\mathcal{M},\mathcal{U})$, will be familiar to those who work with models
of Peano arithmetic: it is the restriction of the usual ultrapower to
equivalence classes of parametrically $\mathcal{M}$-definable functions from
$M$ to $M$. As for ordinary ultrapowers, it is easy to prove a \L o\'{s}
theorem for Skolem ultrapowers: for a function $f$ from $M$ to $M$, a
first-order property holds of the equivalence class $[f]$ of a
parametrically $\mathcal{M}$-definable function in the Skolem ultrapower $%
\mathrm{Ult}(\mathcal{M},\mathcal{U})$ if and only if it holds of $f(i)$ in $%
\mathcal{M}$ for almost all $i$ --- that is, it holds on a set in the
ultrafilter (the definability of Skolem functions is invoked in the
existential case of the inductive proof of the \L o\'{s} theorem to ensure
that the Skolem ultrapower includes the equivalence class of the appropriate
witnessing function). Clearly $\mathrm{Ult}(\mathcal{M},\mathcal{U})$ is
generated over $M$ by the equivalence class of the identity function, $[id]$%
, since for every definable unary function $f$, we have $\mathrm{Ult}(%
\mathcal{M},\mathcal{U})\models [f]=f([id])$ (for example, by applying the
above \L o\'{s} theorem for Skolem ultrapowers). Note that if $\mathcal{U}$\
is an ultrafilter over $\mathcal{P}(M)$, then the ordinary (full) ultrapower
of $\mathcal{M}$ modulo $\mathcal{U}$ coincides with the $\mathcal{L}(%
\mathcal{M})$-reduct of the Skolem ultrapower of $\mathcal{M}^{\#}$ (the
full expansion of $\mathcal{M}$) modulo $\mathcal{U}$.\medskip

We will extend the Skolem ultrapower construction so that instead of a
single generator, we obtain an elementary extension $\mathcal{M}^{\ast }$ of
$\mathcal{M}$ generated by a set of order indiscernibles $(I,<)$, where $I$
is disjoint from $M$. In the degenerate case, where $I$ has a single element
$i$, this \emph{dimensional Skolem ultrapower} is isomorphic to the usual
Skolem ultrapower, where $[i]$ corresponds to $[id]$. Moreover, just as $%
[id] $ forms a single generator of the Skolem ultrapower over $M$, the set
of all $[i]$ (for $i\in I$) forms a set of generators for the dimensional
Skolem ultrapower over $M$. We can identify $[i]$ with $i$ (really, just
renaming). Then we will show that $I$ is an ordered set of tight
indiscernibles. Thus, we will define a notion of \emph{almost all} for
finite sequences from $M$, where the single generator, $[id]$, is replaced
by the generators $[i]$ for $i$ in $I$, so that the following \L o\'{s}
theorem holds: for $i_{1}<...<i_{k}$ from $I$, a first-order property holds
of $\left\langle [i_{1}],\ldots ,[i_{k}]\right\rangle $ in $\mathcal{M}%
^{\ast }$ if and only if for almost all sequences $m_{1}<\ldots <m_{k}$ from
$M^{k}$, the property holds in $\mathcal{M}$. It follows immediately (again,
identifying $[i]$ with $i$, for each $i\in I$) that $I$ is a set of order
indiscernibles in the dimensional Skolem ultrapower.

\subsection{Finite dimensional ultrafilters}

\label{Finite dimensional ultrafilters}

Suppose $\mathcal{U}$ provides a notion of `almost all' as an ultrafilter on
the parametrically $\mathcal{M}$-definable subsets of $M$, and $n$ is a
positive integer. We wish to introduce a notion of `almost all', $\mathcal{U}%
^{n}$, as an ultrafilter on the parametrically $\mathcal{M}$-definable
subsets of $M^{n}.$ Before treating the general case, we warm-up with the
important case of $n=2$. Thus, let $I$ be the two-element order $\{0,1\}$
with $0<1$. Then a parametrically $\mathcal{M}$-definable set $X\subseteq
M^{2}$ belongs to $\mathcal{U}^{2}$ precisely when for almost all $m_{1}\in M
$, it is the case that for almost all $m_{2}\in M$ the pair $\left\langle
m_{1},m_{2}\right\rangle $ is in $X$.\footnote{%
Even though a definable pairing function is available in the archetypical
case of models of arithmetic, our development here does not require that $%
\mathcal{M}$ has a definable pairing function.} To make that precise, it is
convenient to introduce notation for a \emph{section} of $X$, so let:

\begin{center}
$X|m_{1}=$ $\left\{ m_{2}\in M:\left\langle m_{1},m_{2}\right\rangle \in
X\right\} $.
\end{center}

\noindent Then we define $\mathcal{U}^{2}$ as consisting of parametrically $%
\mathcal{M}$-definable subsets of $M^{2}$ such that:

\begin{center}
$\left\{ m_{1}\in M:X|m_{1}\in \mathcal{U}\right\} \in \mathcal{U}$.
\end{center}

\noindent Of course, $X|m_{1}$ is a reasonable candidate for membership in $%
\mathcal{U}$ since $X|m_{1}$ is parametrically $\mathcal{M}$-definable, by
parametric $\mathcal{M}$-definability of $X$. But for the definition above
to yield an ultrafilter, we also need $\left\{ m_{1}\in M:X|m_{1}\in
\mathcal{U}\right\} $ to be parametrically $\mathcal{M}$-definable; else
neither that set nor its complement would be in $\mathcal{U}$. We will need
that definition to be suitably uniform in $m_{1}$ in our formation of $%
\mathcal{U}^{3}$, $\mathcal{U}^{4}$, and so on. We capture this requirement
in the following definition.\medskip

\noindent \textbf{3.1.~Definition.~}Let $\mathcal{M}$ be a structure and $%
\mathcal{U}$ be a non-principal ultrafilter over the parametrically $%
\mathcal{M}$-definable subsets of $M$. $\mathcal{U}$ is $\mathcal{M}$-\emph{%
amenable} if for every first-order $\mathcal{L}(\mathcal{M})$-formula $%
\varphi \left( x,y_{1},\ldots ,y_{k}\right) $, there is a corresponding $%
\mathcal{L}(\mathcal{M})$-formula $U_{\varphi }\left( y_{1},\ldots
,y_{k}\right) $ such that for all $m_{1},\ldots ,m_{k}\in M:$

\begin{center}
$\varphi _{1}^{\mathcal{M}}\in \mathcal{U}$ iff $\mathcal{M}\models
U_{\varphi }(m_{1},\ldots ,m_{k})$,
\end{center}

\noindent where $\varphi _{1}(x)$ denotes $\varphi \left( x,m_{1},\ldots
,m_{k}\right) $. $\mathcal{M}$ is \emph{amenable} if $\mathcal{M}$ has
definable Skolem functions and there is an $\mathcal{M}$-amenable
ultrafilter. \medskip

\noindent \textbf{3.1.1.~Remark.~}Note that if $\mathcal{U}$ is a
nonprincipal ultrafilter over $\mathcal{P}(M)$, and $\mathcal{M}^{\#}$ is a
full expansion of $\mathcal{M}$, then $\mathcal{U}$ is $\mathcal{M}^{\#}$%
-amenable. \medskip

\noindent \textbf{3.2.~Theorem.~}(Amenable Expansions)\textbf{\ }\textit{%
Every infinite structure }$\mathcal{M}$ \textit{has an amenable expansion }$%
\mathcal{M}^{\prime }$\textit{\ with} $\left\vert \mathcal{L}(\mathcal{M}%
^{\prime })\right\vert =\max \left\{ \left\vert \mathcal{L}(\mathcal{M}%
)\right\vert ,\aleph _{0}\right\} .$ \textit{Indeed given any} \textit{%
nonprincipal ultrafilter }$\mathcal{U}\subseteq \mathcal{P}(M)$\textit{\
there is an expansion }$\mathcal{M}^{\prime }$\textit{\ of }$\mathcal{M}$
\textit{with definable Skolem functions such that} $\left\vert \mathcal{L}(%
\mathcal{M}^{\prime })\right\vert =\max \left\{ \left\vert \mathcal{L}(%
\mathcal{M})\right\vert ,\aleph _{0}\right\} $\textit{, and some }$\mathcal{U%
}^{\prime }\subseteq \mathcal{U}$ \textit{such that} $\mathcal{U}^{\prime}$
\textit{is an} $\mathcal{M}^{\prime }$\textit{-amenable ultrafilter}.
\textit{\medskip }

\noindent \textbf{Proof.} Consider the two-sorted structure $\overline{%
\mathcal{M}}=\left( \mathcal{M},\mathcal{P}(M),\pi ,\in ,\mathcal{U}\right) $%
, where $\pi $ is a pairing function on $M$, i.e., a bijection between $M$
and $M^{2}$; $\in $ is the membership relation between elements of $M$ and
elements of $\mathcal{P}(M)$, and $\mathcal{U}$ is construed as a unary
predicate on $\mathcal{P}(M)$. Note that the presence of $\pi $ assures us
that for each positive $k\in \omega $ the structure $\overline{\mathcal{M}}$
carries a definable bijection between $M$ and $M^{k}.$ Since we are assuming
that $\mathrm{ZFC}$ holds in our metatheory, $\overline{\mathcal{M}}$
satisfies the schemes $\Sigma =\left\{ \sigma _{k}:k\in \omega \right\} $
and $A=\left\{ \alpha _{k}:k\in \omega \right\} ,$ where $\sigma _{k}$ is
the $\mathcal{L}(\overline{\mathcal{M}})$-sentence expressing:

\begin{center}
$\forall R\subseteq M^{k+1}\ \exists f:M^{k}\rightarrow M$ $\forall y_{1}\in
M\ldots \forall y_{k}\in M\left( \exists x\ R(x,y_{1},\ldots
,y_{k})\rightarrow R(f(y_{1},\ldots ,y_{k}),y_{1},\ldots ,y_{k})\right) ,$
\end{center}

\noindent and $\alpha _{k}$ is the $\mathcal{L}(\overline{\mathcal{M}})$%
-sentence expressing:

\begin{center}
$\forall R\subseteq M^{k+1}\ \exists U\subseteq M^{k}\ \forall y_{1}\in
M\ldots \forall y_{k}\in M$\quad $\mathcal{U}\left( \{x:R(x,y_{1},\ldots
,y_{k})\}\right) \leftrightarrow \left\langle y_{1},\ldots
,y_{k}\right\rangle \in U.$
\end{center}

\noindent By the L\"{o}wenheim-Skolem theorem there is some $\mathcal{P}%
_{0}(M)\subseteq \mathcal{P}(M)$ of cardinality $\max \left\{ \left\vert%
\mathcal{L}(\mathcal{M})\right\vert ,\aleph _{0}\right\}$ such that:

\begin{center}
$\overset{\overline{\mathcal{M}}_{0}}{\overbrace{\left( \mathcal{M},\mathcal{%
P}_{0}(M),\pi ,\in ,\underset{\mathcal{U}_{0}}{\underbrace{\mathcal{U\cap P}%
_{0}(M)}}\right) }}\prec \overline{\mathcal{M}}.$
\end{center}

\noindent In particular $\overline{\mathcal{M}}_{0}$ satisfies both schemes $%
\Sigma$ and $A,$ and therefore the $\mathcal{M}$-expansion $\mathcal{M}%
^{\prime }=(\mathcal{M},\pi ,X)_{X\in \mathcal{P}_{0}(M)}$ has definable
Skolem functions, and $\mathcal{U}_{0}$ is $\mathcal{M}^{\prime }$%
-amenable.\hfill $\square \medskip $

\noindent \textbf{3.3.~Exercise.~}If $\mathcal{U}$ is an $\mathcal{M}$%
-amenable ultrafilter, then $\mathcal{U}^{2}$ is an ultrafilter on the
parametrically $\mathcal{M}$-definable subsets of $M^{2}$. \smallskip

\noindent \textbf{3.4.~Remark.~}Every model of Peano arithmetic is amenable;
this fact is implicit in the usual proof of the MacDowell-Specker Theorem
\cite[Theorem 2.2.8]{Roman-Jim's book}. Another theory with this property is
$\mathrm{ZF}+\mathrm{V=OD}$, i.e., Zermelo-Fraenkel set theory plus the
axiom expressing that every set is definable from some ordinal. This follows
from three well-known facts: (1) models of $\mathrm{ZF}$ with parametrically
definable Skolem functions are precisely those in which $\mathrm{ZF}+\mathrm{%
V=OD}$ holds, (2) the (global) axiom of choice is provable in $\mathrm{%
ZF+V=OD,}$ and (3) within $\mathrm{ZFC}$ there is a nonprincipal ultrafilter
on any prescribed infinite set. We elaborate (3): Let $\mathcal{M}%
=(M,E)\models \mathrm{ZFC}$, where $E=\ \in ^{\mathcal{M}}$. Fix an infinite
$s\in M$. Then there is some $u\in M$ such that $\mathcal{M}\models $
\textquotedblleft $u$ is a nonprincipal ultrafilter on the power set of $s$%
\textquotedblright . The desired $\mathcal{M}$-amenable ultrafilter $%
\mathcal{U}$\ is the collection of parametrically $\mathcal{M}$-definable
subsets $X$ of $M$ such that $s_{E}\cap X=a_{E}$ for some $a\in u_{E},$
where $x_{E}:=\left\{ y\in M:yEx\right\} .$\medskip

Our next task is to define $\mathcal{U}^{n}$ for arbitrary positive integers
$n$. The lemma following this definition shows that amenability extends
appropriately to powers $n>2$.\medskip

\noindent \textbf{3.5.~Definition.~}For an $\mathcal{M}$-amenable
ultrafilter $\mathcal{U}$, we recursively define $\mathcal{U}^{n}$ for
positive integers $n$: $\mathcal{U}^{1}$ is $\mathcal{U}$, and

\begin{center}
$\mathcal{U}^{n+1}=\{X\subseteq M^{n+1}:X$ is parametrically $\mathcal{M}$%
-definable and $\{m\in M:X|m\in \mathcal{U}^{n}\}\in \mathcal{U}\}.$\medskip
\end{center}

\noindent \textbf{3.6.~Lemma.~}(Extended Amenability)\textbf{\ }\textit{Let}
$\mathcal{U}$ \textit{be an} $\mathcal{M}$\textit{-amenable ultrafilter.
Then for every first-order }$\mathcal{L}(\mathcal{M})$-\textit{formula} $%
\varphi (x_{1},\ldots ,x_{n},y_{1},\ldots ,y_{k})$\textit{, there is a
corresponding }$\mathcal{L}(\mathcal{M})$-\textit{formula }$U_{\varphi
}(y_{1},\ldots ,y_{k})$ \textit{such that for all} $m_{1},\ldots ,m_{k}\in
M: $

\begin{center}
$\varphi _{n}^{\mathcal{M}}\in \mathcal{U}^{n}$ \textit{iff} $\mathcal{M}%
\models U_{\varphi }\left( m_{1},\ldots ,m_{k}\right) $,
\end{center}

\noindent \textit{where }$\varphi _{n}(x_{1},\ldots ,x_{n})$ \textit{denotes}
$\varphi (x_{1},\ldots ,x_{n},m_{1},\ldots ,m_{k}).$\medskip

\noindent \textbf{Proof.} By induction on $n$. The case $n=1$ is just the
definition of $\mathcal{M}$-amenable. Suppose the lemma holds for $n$, and
consider the formula \newline
$\varphi \left( x_{0},x_{1},\ldots ,x_{n},y_{1},\ldots ,y_{k}\right) $. By
the inductive hypothesis, there is a formula $U_{\varphi
}(x_{0},y_{1},\ldots ,y_{k})$ such that for all $m,m_{1},\ldots ,m_{k}\in M$%
, we have:

\begin{center}
$\varphi _{m}^{\mathcal{M}}\in \mathcal{U}^{n}$ iff $\mathcal{M}\models
U_{\varphi }(m,m_{1},\ldots ,m_{k})$,
\end{center}

\noindent where $\varphi_{m}=\varphi (m,x_{1},...,x_{n},m_{1},...,m_{k})$.
By amenability, there is a formula $V_{\varphi }(y_{1},\ldots ,y_{k})$ such
that for all $m_{1},\ldots ,m_{k}\in M$ and for $\varphi
_{0}(x_{0})=U_{\varphi }(x_{0},m_{1},\ldots ,m_{k}):$

\begin{center}
$\varphi_{0}^{\mathcal{M}}\in \mathcal{U}$ iff $\mathcal{M}\models
V_{\varphi }(m_{1},\ldots,m_{k})$.
\end{center}

\noindent Then setting $\varphi ^{+}(x_{0},x_{1},\ldots,x_{n})=\varphi
(x_{0},x_{1},\ldots,x_{n},m_{1},\ldots,m_{k})$, we have the following
equivalences:

\begin{center}
$\left( \varphi ^{+}\right) ^{\mathcal{M}}\in \mathcal{U}^{n+1}$ \smallskip

$\Updownarrow $ (by the definition of $\mathcal{U}^{n+1},$ $\varphi ^{+},$
and $\varphi _{m})$ \smallskip

$\left\{ m\in M:\varphi _{m}^{\mathcal{M}}\in \mathcal{U}^{n}\right\} \in
\mathcal{U}$ \smallskip

$\Updownarrow $ (by the above choice of $U_{\varphi })$\smallskip

$\left\{ m\in M:\mathcal{M}\models U_{\varphi
}(m,m_{1},\ldots,m_{k})\right\} \in \mathcal{U}$ \smallskip

$\Updownarrow $\ (by the above choice of $V_{\varphi })$\smallskip

$\mathcal{M}\models V_{\varphi }(m_{1},\ldots ,m_{k}).$

\hfill $\square \medskip $
\end{center}

\noindent \textbf{3.7.~Exercise~}(Finite dimensions provide ultrafilters).
Let $\mathcal{U}$ be an $\mathcal{M}$-amenable ultrafilter on the
parametrically $\mathcal{M}$-definable subsets of $M$. Then for all positive
integers $n$, $\mathcal{U}^{n}$ is an ultrafilter on the parametrically $%
\mathcal{M}$-definable subsets of $M^{n}$.\smallskip

\subsection{Warm-up: two-dimensional Skolem ultrapowers}

\label{Warmup on two dim. ultrapowers}

\begin{itemize}
\item We assume throughout the rest of this section that $\mathcal{M}$ has
definable Skolem functions, and $\mathcal{U}$ is an $\mathcal{M}$-amenable
ultrafilter.
\end{itemize}

Now that we have defined $\mathcal{U}^{n}$ for $\mathcal{M}$-amenable
ultrafilters $\mathcal{U}$, we want to use them to construct the desired
extension of $\mathcal{M}$ by tight indiscernibles $(I,<)$. Let us begin by
taking a close look at the simple but important case of $n=2.$ Consider the
ordered set $I=2=\{0,1\}$, where $0<1$. To form the 2-dimensional Skolem
ultrapower $\mathcal{M}^{\ast }=\mathrm{Ult}(\mathcal{M},\mathcal{U},2)$,
instead of taking equivalence classes of \textit{unary} functions as we
would when constructing an ordinary Skolem ultrapower, we take equivalence
classes of parametrically $\mathcal{M}$-definable \textit{binary} functions,
where the equivalence relation at work is :

\begin{center}
$f(x,y)\thicksim g(x,y)$ iff $\left\{ \left\langle x,y\right\rangle
:f(x,y)=g(x,y)\right\} \in \mathcal{U}^{2}.$
\end{center}

\noindent We write $[f(0,1)]$ to denote the equivalence class of $f$, where
here $0$ and $1$ refer to elements of the order, $\{0,1\}$. The universe of
the desired dimensional ultrapower $\mathcal{M}^{\ast }$ is the set of all
such $[f(0,1)]$. Two distinguished examples of such objects are the case of $%
f=id_{0}$ and $f=id_{1}$ where $id_{0}(x,y)=x$ and $id_{1}(x,y)=y$. For
convenience, we may write $[0]$ and $[1]$ for these respective equivalence
classes. One can check that $[f(0,1)]=f^{\mathcal{M}^{\ast }}([0],[1])$;
hence, the set $\{[0],[1]\}$ is a set of generators for $\mathcal{M}^{\ast }$%
. But before we can do that, we need to define functions and relations on $%
\mathcal{M}^{\ast }$, i.e., on this set of equivalence classes. For each
function symbol $g$ and arguments $[f_{i}(0,1)]$ of $\mathcal{M}^{\ast }$
whose length is the arity of $g$, interpret $g$ in $\mathcal{M}^{\ast }$ on
these arguments by applying $g$ pointwise in $\mathcal{M}$:

\begin{center}
$g([f_{0}(0,1)],...,[f_{k-1}(0,1)])=[h(0,1)],$
\end{center}

\noindent where $h(x,y)=g(f_{0}(x,y),...,f_{k-1}(x,y))$. This is
well-defined: $h$ is parametrically $\mathcal{M}$-definable, and if $%
[f_{i}(0,1)]=[f_{i}^{\prime }(0,1)]$ for $i<k$, then $f_{i}(x,y)=f_{i}^{%
\prime}(x,y)$ on a set in the ultrafilter $\mathcal{U}^{2}$ for each such $i$%
, hence for all $i$ in some $X\in \mathcal{U}^{2}$; and on this set $X$ of
pairs, $h(x,y)$ is unchanged if we replace each $f_{i}$ by $f_{i}^{\prime}$
in the definition of $h$. \medskip

\noindent \textbf{3.8.~Exercise~(a) }For the case of $\mathcal{U}^{2}$
above, complete the definition of $\mathcal{M}^{\ast }$ (by interpreting
relation symbols) and then prove the \L o\'{s} theorem. \textbf{(b)} Use
part (a) to show that $\left\{ [0],[1]\right\} $ forms a set of order
indiscernibles in $\mathcal{M}^{\ast }$ over $\mathcal{M}$.

\subsection{Dimensional Skolem ultrapowers}

\label{Dimensional Ultrapowers}

Assume $(I,<)$ is an ordered set disjoint from $M$, and $\mathcal{U}$ is an $%
\mathcal{M}$-amenable ultrafilter. In order to define $\mathcal{M}^{\ast }=%
\mathrm{Ult}(\mathcal{M},\mathcal{U},I)$ (the $I$-dimensional Skolem
ultrapower of $\mathcal{M}$ with respect to $\mathcal{U)}$ we extend the
notion of dimensional ultrapower from $I=\{0,1\}$ to arbitrary ordered sets $%
I$. For all non-empty finite subsets $I_{0}$ of $I$, we call the function
space $I_{0}\rightarrow M$, also denoted $M^{I_{0}}$, the \emph{set of $%
I_{0} $-sequences }(\emph{from }$M)$. Given a finite subset $I_{0}$ of $I$,
where $\left\langle i_{0},\ldots ,i_{k}\right\rangle $ enumerates $I_{0}$ in
increasing order, we can take advantage of the order isomorphism $%
j$ from $\{0,...,k\}$ to $I_{0}$ to form a `copy' $\mathcal{U}^{I_{0}}$ of $%
\mathcal{U}^{\left\vert I_{0}\right\vert }$ on the set of $I_{0}$-sequences,
where as usual $\left\vert \ \right\vert $ denotes cardinality. More
specifically, for $X\subseteq M^{I_{0}}$, let $X[j]\subseteq M^{k+1}$ be $%
\{s\circ j:s\in X\}$, and then define:

\begin{center}
$\mathcal{U}^{I_{0}}=\left\{ X\subseteq M^{I_{0}}:X[j]\in \mathcal{U}%
^{k+1}\right\}$.
\end{center}

\noindent For finite $I_{0}\subseteq I$ and function $f$ with arity $%
\left\vert I_{0}\right\vert $, we write $f(I_{0})$ as an abbreviation for
the \textit{syntactic entity}

\begin{center}
$f(i_{0},\ldots ,i_{k}),$
\end{center}

\noindent where $f$ is a $k+1$-ary parametrically $\mathcal{M}$-definable
function and $\left\langle i_{0},\ldots ,i_{k}\right\rangle $ enumerates $%
I_{0}$ in increasing order; we call this a \emph{generalized term} (over $%
I_{0}$). Moreover, if $u$ is an $I_{1}$-sequence for some $I_{1}$ containing
$I_{0}$, we write $f\left( i_{0},\ldots ,i_{k}\right) [u]$ to denote the
value $f^{\mathcal{M}}(u(i_{0}),\ldots ,u(i_{k}))$, where $f^{\mathcal{M}}$
is the interpretation of $f$ in $\mathcal{M}$. Given finite $I_{0}$, $%
I_{1},I_{2}\subseteq I$, with $I_{2}=I_{0}\cup I_{1}$, and generalized terms
$f(I_{0})$ and $g(I_{1})$, where $f$ and $g$ are allowed to have different
arities, we say that $f(I_{0})$ and $g(I_{1})$ are $\mathcal{U}^{\ast }$-%
\emph{equivalent} if

\begin{center}
$\left\{ u\in M^{I_{2}}:f(I_{0})[u]=g(I_{1})[u]\right\} \in \mathcal{U}%
^{I_{2}}$.
\end{center}

\noindent The universe $M^{\ast }$ of the $I$-dimensional ultrapower $%
\mathcal{M}^{\ast }$ with respect to $\mathcal{U}$ is the set of all $%
\mathcal{U}^{\ast }$-equivalence classes $\left[ f(I_{0})\right] $, where $%
I_{0}$ $\subseteq I$, $\left\vert I_{0}\right\vert =n\in \omega $, and $f$
is an $n$-ary parametrically $\mathcal{M}$-definable function. Note that $M$
naturally can be injected into $M^{\ast }$ by:

\begin{center}
$\varepsilon :M\rightarrow M^{\ast }$,
\end{center}

\noindent where $\varepsilon (m)$ is defined as $[f_{m}(I_{0})]$, $f_{m}$ is
the constant function whose range is $\{m\}$, and $I_{0}$ is any finite
subset of $I.$ Functions and relations interpreted in $\mathcal{M}$ extend
naturally to $\mathcal{M}^{\ast }$: for each $n+1$-ary function symbol $g$
of $\mathcal{L}(\mathcal{M})$, by defining

\begin{center}
$g^{\mathcal{M}^{\ast }}([f_{0}(I_{0})],[f_{1}(I_{1})],\ldots
,[f_{n}(I_{n})])=[h(I_{n+1})],$
\end{center}

\noindent where $I_{0},...,I_{n}$ are allowed to overlap, $%
I_{n+1}=\bigcup\limits_{0\leq j\leq n}I_{j}$, and $h$ is a parametrically $%
\mathcal{M}$-definable function such that for all $u\in M^{I_{n+1}}$,

\begin{center}
$h(I_{n+1})[u]=g^{\mathcal{M}}(f_{0}(I_{0})[u],\ldots ,f_{k}(I_{n})[u])$.
\end{center}

\noindent \textbf{3.9.~Exercise.}~Define the interpretation of a relation
symbol in $\mathcal{M}^{\ast }$ in analogy to how function symbols are
interpreted (as above). Then verify that $\varepsilon $ is an isomorphism
between $\mathcal{M}$ and $\varepsilon (\mathcal{M}),$ where $\varepsilon (%
\mathcal{M})$ is the submodel of $\mathcal{M}^{\ast }$ whose universe is $%
\varepsilon (M).$\medskip

The isomorphism of $\mathcal{M}$ with $\varepsilon (\mathcal{M})$ allows
us---as is commonly done in model theory---to replace $\mathcal{M}^{\ast }$
with an isomorphic copy so as to arrange $\mathcal{M}$ to be \textit{submodel%
} of $\mathcal{M}^{\ast }$. Therefore we may identify any element of $%
\mathcal{M}^{\ast }$ of the form $[f_{m}(I_{0})]$ with $m$ itself. With $%
\mathcal{M}^{\ast }$ defined, many of its properties now fall naturally into
place, as explained below.\smallskip

It is useful to extend the notation $X|x$ to $X|s$ for an ordered sequence $%
s $. Assume that $I_{0}$, $I_{1}\subseteq I$ such that $\max (I_{0})<\min
(I_{1})$; then let $I_{2}=I_{0}\cup I_{1}$, and assume that $X\subseteq
M^{I_{2}}$ and $s\in M^{I_{0}}$. We want to define a subset $X|s$ of $%
M^{I_{1}}$ to be the result of collecting, for each sequence in $X$ that
starts with $s$, its restriction to $I_{1}$. The recursive definition of $%
X|s $ is as follows, where $\left\langle \ \right\rangle $ is the empty
sequence. $X|\left\langle \ \right\rangle $ is $X$. Now suppose $I_{0}$ is
ordered as $i_{0}<\ldots <i_{k}$ and let $I_{0}^{\prime }=\{i_{1},\ldots
,i_{k}\}$. Then $X|I_{0}$ is $(X|i_{0})|I_{0}^{\prime}$. The following
exercise generalizes the recursive definition of $\mathcal{U}^{n}$ by
stripping off initial subsequences rather than merely single initial
elements.\medskip

\noindent \textbf{3.10.~Exercise.}~Let $I_{0}$ be a finite subset of $I$,
and suppose $I_{0}=I_{1}\cup I_{2}$, where $\max (I_{1})<\min (I_{2}$). Then:

\begin{center}
$X\in \mathcal{U}^{I_{0}}$ iff $\left\{ s\in M^{I_{1}}:X|s\in \mathcal{U}%
^{I_{2}}\right\} \in \mathcal{U}^{I_{1}}$.
\end{center}

\noindent We next establish a useful lemma.\medskip

\noindent \textbf{3.11.~Lemma.}~\textit{Suppose} $I_{0}$ \textit{is a finite
subset of} $I$, \textit{and} $J$ \textit{is a finite subset of} $I$ \textit{%
with} $I_{0}\subseteq J.$ \textit{If} $X\subseteq M^{I_{0}}$ \textit{and}

\begin{center}
$X^{\prime }=\left\{ s\cup t\in M^{J}:s\in X\wedge t\in M^{J\backslash
I_{0}}\right\},$
\end{center}

\noindent \textit{then}

\begin{center}
$X\in \mathcal{U}^{I_{0}}$ iff $X^{^{\prime }}\in \mathcal{U}^{J}.$
\end{center}

\noindent \textbf{Proof.} We prove the Lemma when $|J\backslash I_{0}|=1$.
The Lemma in full generality then follows by induction on $|J\backslash
I_{0}|$. So, suppose $J\backslash I_{0}=\{i\}$. Let $X\subseteq M^{I_{0}}$,
and let
\begin{equation*}
X^{\prime }=\{s\cup t\in M^{J}:s\in X\wedge t\in M^{J\backslash I_{0}}\}.
\end{equation*}%
There are three cases to consider:\newline
\textbf{Case} $i<\min I_{0}$: By Exercise 3.10,
\begin{equation*}
X^{\prime }\in \mathcal{U}^{J}\text{ if and only if }\{s\in
M^{\{i\}}:X^{\prime }|s\in \mathcal{U}^{I_{0}}\}\in \mathcal{U}^{\{i\}}.
\end{equation*}%
Now, for all $s\in M^{\{i\}}$, $X^{\prime }|s=X$, so
\begin{equation*}
\{s\in M^{\{i\}}:X^{\prime }|s\in \mathcal{U}^{I_{0}}\}=\left\{
\begin{array}{ll}
M^{\{i\}} & \text{if }X\in \mathcal{U}^{I_{0}} \\
\emptyset & \text{otherwise}%
\end{array}%
\right.
\end{equation*}%
Therefore, $X^{\prime }\in \mathcal{U}^{J}$ if and only if $X\in \mathcal{U}%
^{I_{0}}$.\newline
\textbf{Case} $\max I_{0}<i$: By Exercise 3.10,
\begin{equation*}
X^{\prime }\in \mathcal{U}^{J}\text{ if and only if }\{s\in
M^{I_{0}}:X^{\prime }|s\in \mathcal{U}^{\{i\}}\}\in \mathcal{U}^{I_{0}}.
\end{equation*}%
Now, for all $s\in M^{I_{0}}$,
\begin{equation*}
X^{\prime }|s=M^{\{i\}}\text{ if and only if }s\in X.
\end{equation*}%
Therefore, $X^{\prime }\in \mathcal{U}^{J}$ if and only if $X\in \mathcal{U}%
^{I_{0}}$.\newline
\textbf{Case} $(\exists i^{\prime },i^{\prime \prime }\in I_{0})(i^{\prime
}<i<i^{\prime \prime })$: Let $I_{0}=I_{1}\cup I_{2}$ with $\max
I_{1}<i<\min I_{2}$. By Exercise 3.10,
\begin{equation*}
X^{\prime }\in \mathcal{U}^{J}\text{ if and only if }\{s\in
M^{I_{1}}:X^{\prime }|s\in \mathcal{U}^{\{i\}\cup I_{2}}\}\in \mathcal{U}%
^{I_{1}}.
\end{equation*}%
By the first case that we considered, for all $s\in M^{I_{1}}$,
\begin{equation*}
X^{\prime }|s\in \mathcal{U}^{\{i\}\cup I_{2}}\text{ if and only if }X|s\in
\mathcal{U}^{I_{2}}.
\end{equation*}%
Therefore, using Exercise 3.10,
\begin{equation*}
X^{\prime }\in \mathcal{U}^{J}\text{ if and only if }\{s\in M^{I_{1}}:X|s\in
\mathcal{U}^{I_{2}}\}\in \mathcal{U}^{I_{1}}
\end{equation*}%
\begin{equation*}
\text{if and only if }X\in \mathcal{U}^{I_{0}}.
\end{equation*}%
\hfill $\square \medskip $

\noindent \textbf{3.12.~Exercise.}~Use Lemma 3.11 to show that the
definition of $\mathcal{U}^{\ast }$-equivalence above is unchanged if $I_{2}$
is replaced by any finite superset of $I_{2}$ (i.e., of $I_{0}\cup I_{1}$)
contained in $I$.\smallskip

\noindent \textbf{3.13.~Exercise.}~Formulate and prove the \L o\'{s} theorem
for $\mathcal{M}^{\ast }$, and use it to verify that $\varepsilon $ is an
elementary embedding. \smallskip

\noindent \textbf{3.14.~Exercise.}~Prove property (1) of tight
indiscernibles for $\mathcal{M}^{\ast }$ and $(I,<)$.\smallskip

It remains to prove properties (2) and (3) of tight indiscernibles for $%
\mathcal{M}^{\ast }$ and $(I,<)$. Property (2) is easy to see since if $%
[f(i_{0},\ldots ,i_{k})]\in M^{\ast },$ then by the very manner in which
functions symbols of $\mathcal{L}(\mathcal{M})$ are interpreted in $\mathcal{%
M}^{\ast }$, we have:

\begin{center}
$[f(i_{0},\ldots ,i_{k})]=f^{\mathcal{M}^{\ast }}(i_{0},\ldots ,i_{k}).$
\end{center}

\noindent Finally, to verify property (3) of tight indiscernibles for $%
\mathcal{M}^{\ast}$ and $(I,<)$ suppose $I_{0}$ and $I_{1}$ are distinct
finite subsets of $I$ where $\max (I_{0})<$ $\min (I_{1})$ of $I_{1}$, and
suppose that $f(I_{0})$ and $f(I_{1})$ are generalized terms such that $%
\left[ f(I_{0})\right] =\left[ f(I_{1})\right] $; we show that this common
value is in $M$. Let us first consider the simple case that $I_{0}=\{i_{0}\}$
and $I_{1}=\{i_{1}\}$, with $i_{0}<i_{1}$. Let

\begin{center}
$X=\left\{ \left\langle x,y\right\rangle \in M^{2}:f^{\mathcal{M}}(x)=f^{%
\mathcal{M}}(y)\right\} $;
\end{center}

\noindent then from $[f(i_{0})]=[f(i_{1})]$ and the definition of $\mathcal{U%
}^{\ast }$-equivalence, $X\in \mathcal{U}^{\{i_{0},i_{1}\}}$; thus,

\begin{center}
$\{x_{0}\in M:X|x_{0}\in \mathcal{U}\}\in \mathcal{U}$.
\end{center}

\noindent In particular, we may choose $x_{0}\in M$ such that $X|x_{0}\in
\mathcal{U}$. Let $r=f^{\mathcal{M}}(x_{0})$; then $r\in M$ and for all $%
y\in X|x_{0}$, $f(y)=r$. Therefore by definition, $[f(I_{1})]\in M$, so $%
m\in M$. \medskip

Enumerate $I_{0}$ as $i_{1}<\ldots <i_{k}$ and enumerate $I_{1}$ as $%
j_{1}<\ldots <j_{k}$, where $I_{0},I_{1}\subseteq I$ and $i_{k}<j_{1}$, and
assume that $f$ is a parametrically $\mathcal{M}$-definable function such
that $f(i_{1},\ldots ,i_{k})=f(j_{1},\ldots ,j_{k}).$ Let $I_{2}=$ $%
I_{0}\cup I_{1}$ and let

\begin{center}
$X=\left\{ \left\langle x_{1},\ldots,x_{k},y_{1},\ldots,y_{k}\right\rangle
\in M^{I_{2}}:f(x_{1},\ldots,x_{k})=f(y_{1},\ldots,y_{k})\right\} $;
\end{center}

\noindent then by definition of $\mathcal{U}^{\ast }$-equivalence and the
equality of $f(i_{1},\ldots ,i_{k})$ and $f(j_{1},\ldots ,j_{k})$, $X\in
\mathcal{U}^{I_{2}}$. By Exercise 3.10 and the fact that each element of $%
\mathcal{U}^{I_{0}}$ is non-empty, we may choose $s\in M^{I_{0}}$ such that $%
X|s\in \mathcal{U}^{I_{1}}$. Let $r\in M$ such that $r=f(s)$. Thus

\begin{center}
$\left\{ \left\langle y_{1},\ldots,y_{k}\right\rangle \in
M^{I_{1}}:f(y_{1},\ldots,y_{k})=r\right\} \in \mathcal{U}^{I_{1}} $,
\end{center}

\noindent hence $[f(j_{1},\ldots ,j_{k})]=r\in M$. \hfill $\square $

\section{History and other applications}

\label{History}

Skolem \cite{Skolem} introduced the definable ultrapower construction to
exhibit a nonstandard model of arithmetic. Full ultraproducts/powers were
invented by \L o\'{s} in his seminal paper \cite{Jerzy}. The history of
iterated ultrapowers is concisely summarized by Chang \& Keisler in their
canonical textbook \cite[Historical Notes for 6.5]{Chang-Keisler} as follows:

\begin{quote}
{\small Finite iterations of ultrapowers were developed by Frayne, Morel,
and Scott \cite{Frayne et al.}. The infinite iterations were introduced by
Gaifman \cite{Haim-Unform}. Our presentation is a simplification of
Gaifman's work. Gaifman used a category theoretic approach instead of a
function which lives on a finite set. Independently, Kunen \cite{Ken-L(U)}
developed iterated ultrapowers in essentially the same way as in this
section, and generalized the construction even further to study models of
set theory and measurable cardinals.}
\end{quote}

A thorough treatment of (full) iterated ultrapowers, including the many
important results that are relegated to the exercises can be found in \cite[%
Sec.~6.5]{Chang-Keisler}; e.g., \cite[Exercise 6.5.27]{Chang-Keisler} asks
the reader to verify that indiscernible sequences arise from the iterated
ultrapower construction. In the general treatment of the subject, one is
allowed to use a different ultrafilter (for the formation of the ultrapower)
at different stages of the iteration process, but we decided to focus our
attention in this paper on the conceptually simpler case, where the same
ultrafilter is used at all stages of the iteration process. It should be
noted, however, that the framework developed in this paper naturally lends
itself to a wider framework to handle iterations using different
ultrafilters at each stage of the iteration.\medskip

Our treatment of iterated (Skolem) ultrapowers over amenable structures has
two sources of inspiration: Kunen's aforementioned adaptation of iterated
ultrapowers for models of set theory (where the ultrafilters at work are
referred to as $\mathcal{M}$-ultrafilters\footnote{%
This terminology is not standard: Kunen's $\mathcal{M}$-ultrafilters are
also referred to as \textit{iterable} ultrafilters or as (\textit{weakly})%
\textit{\ amenable} ultrafilters by other authors.}), and Gaifman's
adaptation \cite{Haim-PA} of iterated ultrapowers for models of arithmetic.
Gaifman's adaptation is recognized as a major tool in the model theory of
Peano arithmetic; see \cite{Roman-Jim's book}, \cite{Ali-Tehran}, \cite%
{Ali-FM} and \cite{Ali-APAL} for applications and generalizations in this
direction. The referee has also kindly reminded us of the fact that
Gaifman's work has intimate links -- at the conceptual, as well as the
historical level -- with the grand subdiscipline of model theory known as
stability theory. For example, Definition 3.1 can be recast in the jargon of
stability theory to assert that an ultrafilter $\mathcal{U}$ on the
parametrically $\mathcal{M}$-definable sets is amenable precisely when $%
\mathcal{U}$ is a `definable type' when viewed as a 1-type.\footnote{%
The referee has also commented that Shelah honed in on the property that
every type is definable, and showed that a complete theory $T$ is stable iff
given any $\mathcal{M}\models T$, \textit{every} ultrafilter over the
parametrically $\mathcal{M}$-definable sets is amenable when viewed as a
1-type. However, it does not follow that every stable structure is amenable,
as adding Skolem functions can destroy stability.} According to Poizat \cite[%
p.237]{Bruno}:

\begin{quote}
{\small Some trustworthy witnesses assert that the notion of definable types
was not introduced in 1968 by Shelah, but by Haim Gaifman, in order to
construct end extensions of models of arithmetic, see \cite{Haim-PA}.}
\end{quote}

\noindent Furthermore, in his review of Gaifman's paper \cite{Haim-PA} ,
Ressayre \cite{Ressayre review} writes:

\begin{quote}
{\small There is a paradoxical link [between Gaifman's paper] with stable
first-order theories: although the notion of definable type was introduced
by Gaifman in the study of PA, which is the most unstable theory, this
notion turned out to be a fundamental one for stable theories (...) I
expect\ (i) that it will not be possible to \textquotedblleft
explain\textquotedblright\ this similarity by a (reasonable) common
mathematical theory; and (ii) that this similarity is not superficial.
Although they cannot be captured mathematically, such similarities do occur
repeatedly and not by chance in the development of two opposite parts of
logic, namely model theory in the algebraic style on one hand, and the
theory (model, proof, recursion, and set theory) of the basic universes
(e.g. arithmetic, analysis, V , etc.) and their axiomatic systems on the
other.}
\end{quote}

\noindent See also \cite[Remark 0.1]{Haim-PA}, describing the relationship
between Gaifman's notion of `end extension type' with the Harnik-Ressayre
notion of `minimal type'{\small .}\medskip

Iterated ultrapowers have proved indispensable in the study of inner models
of large cardinals ever since Kunen's work \cite{Ken-L(U)}; see Jech's
monograph \cite[Chapter 19]{Tomas} for the basic applications, and Steel's
lecture notes \cite{Steel} for the state of the art. The applications of
iterated ultrapowers in set theory are focused on well-founded models of set
theory, however iterated ultrapowers are also an effective tool in the study
of automorphisms of models of set theory (where well-founded models are of
little interest since no well-founded model of the axiom of extensionality
has a nontrivial automorphism); see, e.g., \cite{Mahlo} and \cite%
{Zach+Matt+Ali}. Iterated ultrapowers have also found many applications
elsewhere, e.g., in the work of Hrb\'{a}\v{c}ek \cite{Karel}, and Kanovei \&
Shelah \cite{Kanovei-Shelah} on the foundations of nonstandard analysis.
More recently, iterated ultrapowers have found new applications to the model
theory of infinitary logic by Baldwin and Larson \cite{Baldwin-Larson}, and
by Boney \cite{Boney}. \medskip

\bibliographystyle{alpha}
\bibliography{.}      

\appendix

\section[Solutions to the exercises]{Solutions to the exercises}

\noindent \textbf{3.3. Exercise.} Let $\mathcal{U}$ be an $\mathcal{M}$%
-amenable ultrafilter on the parametrically $\mathcal{M}$-definable subsets
of $M$. We begin by showing that $\mathcal{U}^{2}$ is a filter on the
parametrically $\mathcal{M}$-definable subsets of $M^{2}$. Let $X,Y\subseteq
M^{2}$ be parametrically $\mathcal{M}$-definable. Let
\begin{equation*}
Z=\{m\in M:(X\cap Y)|m\in \mathcal{U}\},
\end{equation*}%
\begin{equation*}
Z_{1}=\{m\in M:X|m\in \mathcal{U}\},
\end{equation*}%
\begin{equation*}
Z_{2}=\{m\in M:Y|m\in \mathcal{U}\}.
\end{equation*}%
Since $\mathcal{U}$ is $\mathcal{M}$-amenable, $Z$, $Z_{1}$ and $Z_{2}$ are
parametrically $\mathcal{M}$-definable subsets of $M$. Now,
\begin{equation*}
Z=\{m\in M:X|m\cap Y|m\in \mathcal{U}\}.
\end{equation*}%
And, since $\mathcal{U}$ is a filter,
\begin{equation*}
X|m\cap Y|m\in \mathcal{U}\text{ if and only if }X|m\in \mathcal{U}\text{
and }Y|m\in \mathcal{U}.
\end{equation*}%
Therefore, using the fact that $\mathcal{U}$ is a filter,
\begin{equation*}
X\cap Y\in \mathcal{U}^{2}\iff Z\in \mathcal{U}\iff Z_{1}\cap Z_{2}\in
\mathcal{U}\iff
\end{equation*}%
\begin{equation*}
Z_{1}\in \mathcal{U}\text{ and }Z_{2}\in \mathcal{U}\iff X\in \mathcal{U}^{2}%
\text{ and }Y\in \mathcal{U}^{2}.
\end{equation*}%
It follows that if $X,Y\in \mathcal{U}^{2}$, then $X\cap Y\in \mathcal{U}%
^{2} $; and if $X\in \mathcal{U}^{2}$ and $X\subseteq Y$, then $Y\in
\mathcal{U}^{2}$. Therefore $\mathcal{U}^{2}$ is a filter. We are left to
verify that $\mathcal{U}^{2}$ is an ultrafilter. Let $X$ and $Z_{1}$ be as
above. Let
\begin{equation*}
W=\{m\in M:(M^{2}\backslash X)|m\in \mathcal{U}\}.
\end{equation*}%
Again, since $\mathcal{U}$ is $\mathcal{M}$-amenable, $W$ is a
parametrically $\mathcal{M}$-definable subset of $M$. Note that
\begin{equation*}
W=\{m\in M:M\backslash (X|m)\in \mathcal{U}\}.
\end{equation*}%
And, since $\mathcal{U}$ is an ultrafilter, for all $m\in M$,
\begin{equation*}
M\backslash (X|m)\in \mathcal{U}\text{ if and only if }X|m\notin \mathcal{U}.
\end{equation*}%
Therefore
\begin{equation*}
W=\{m\in M:X|m\notin \mathcal{U}\}=M\backslash Z_{1}.
\end{equation*}%
And so, since $\mathcal{U}$ is an ultrafilter,
\begin{equation*}
M^{2}\backslash X\in \mathcal{U}^{2}\iff W\in \mathcal{U}\iff Z_{1}\notin
\mathcal{U}\iff X\notin \mathcal{U}^{2}.
\end{equation*}%
This shows that $\mathcal{U}^{2}$ is an ultrafilter on the parametrically $%
\mathcal{M}$-definable subsets of $M^{2}$.

\noindent \textbf{3.7. Exercise.} Let $\mathcal{U}$ be an $\mathcal{M}$%
-amenable ultrafilter on the parametrically $\mathcal{M}$-definable subsets
of $M$. We will prove that for all positive integers $n$, $\mathcal{U}^{n}$
is an ultrafilter on the parametrically $\mathcal{M}$-definable subsets of $%
M^{n}$ by induction on $n$. The proof is an obvious generalization of
Exercise 3.3. Let $n>0$ and suppose that $\mathcal{U}^{n}$ is an ultrafilter
on the parametrically $\mathcal{M}$-definable subsets of $M^{n}$. We need to
show that $\mathcal{U}^{n+1}$ is an ultrafilter on the parametrically $%
\mathcal{M}$-definable subsets of $M^{n+1}$. We begin by showing that $%
\mathcal{U}^{n+1}$ is a filter. Let $X,Y\subseteq M^{n+1}$ be parametrically
$\mathcal{M}$-definable. Let
\begin{equation*}
Z=\{m\in M:(X\cap Y)|m\in \mathcal{U}^{n}\},
\end{equation*}%
\begin{equation*}
Z_{1}=\{m\in M:X|m\in \mathcal{U}^{n}\},
\end{equation*}%
\begin{equation*}
Z_{2}=\{m\in M:Y|m\in \mathcal{U}^{n}\}.
\end{equation*}%
It follows from Lemma 3.6 that $Z$, $Z_{1}$ and $Z_{2}$ are parametrically $%
\mathcal{M}$-definable subsets of $M$. Now,
\begin{equation*}
Z=\{m\in M:X|m\cap Y|m\in \mathcal{U}^{n}\}.
\end{equation*}%
By the induction hypothesis, for all $m\in M$,
\begin{equation*}
X|m\cap Y|m\in \mathcal{U}^{n}\text{ if and only if }X|m\in \mathcal{U}^{n}%
\text{ and }Y|m\in \mathcal{U}^{n}.
\end{equation*}%
Therefore, using the fact that $\mathcal{U}$ is a filter,
\begin{equation*}
X\cap Y\in \mathcal{U}^{n+1}\iff Z\in \mathcal{U}\iff Z_{1}\cap Z_{2}\in
\mathcal{U}\iff
\end{equation*}%
\begin{equation*}
Z_{1}\in \mathcal{U}\text{ and }Z_{2}\in \mathcal{U}\iff X\in \mathcal{U}%
^{n+1}\text{ and }Y\in \mathcal{U}^{n+1}.
\end{equation*}%
It follows that if $X,Y\in \mathcal{U}^{n+1}$, then $X\cap Y\in \mathcal{U}%
^{n+1}$; and if $X\in \mathcal{U}^{n+1}$ and $X\subseteq Y$, then $Y\in
\mathcal{U}^{n+1}$. This shows that $\mathcal{U}^{n+1}$ is a filter. We are
left to verify that $\mathcal{U}^{n+1}$ is an ultrafilter. Let $X$ and $%
Z_{1} $ be as above. Let
\begin{equation*}
W=\{m\in M:(M^{n+1}\backslash X)|m\in \mathcal{U}^{n+1}\}.
\end{equation*}%
It follows from Lemma 3.6 that $W$ is a parametrically $\mathcal{M}$%
-definable subset of $M$. Note that
\begin{equation*}
W=\{m\in M:M^{n}\backslash (X|m)\in \mathcal{U}^{n}\}.
\end{equation*}%
By the induction hypothesis we have, for all $m\in M$,
\begin{equation*}
M^{n}\backslash (X|m)\in \mathcal{U}^{n}\text{ if and only if }X|m\notin
\mathcal{U}^{n}.
\end{equation*}%
Therefore
\begin{equation*}
W=\{m\in M:(X|m)\notin \mathcal{U}^{n}\}=M\backslash Z_{1}.
\end{equation*}%
Using the fact that $\mathcal{U}$ is an ultrafilter, this shows that
\begin{equation*}
M^{n+1}\backslash X\in \mathcal{U}^{n+1}\iff W\in \mathcal{U}\iff
Z_{1}\notin \mathcal{U}\iff X\notin \mathcal{U}^{n+1},
\end{equation*}%
which proves that $\mathcal{U}^{n+1}$ is an ultrafilter on the
parametrically $\mathcal{M}$-definable subsets of $M^{n+1}$.

\noindent \textbf{3.8. Exercise. (a).} Let $\mathcal{U}$ be an $\mathcal{M}$%
-amenable ultrafilter. Let $\mathcal{M}^{\ast }=\mathrm{Ult}(\mathcal{M},%
\mathcal{U},2)$ (defined in Subsection \ref{Warmup on two dim. ultrapowers}%
). We need to define the interpretations of the relation symbols in $%
\mathcal{L}(\mathcal{M})$ in $\mathcal{M}^{\ast }$. Let $R(x_{1},\ldots
,x_{n})$ be an $n$-ary relation symbol in $\mathcal{L}(\mathcal{M})$. For
all $[f_{1}(0,1)],\ldots ,[f_{n}(0,1)]\in M^{\ast }$, define
\begin{equation*}
R^{\mathcal{M}^{\ast }}([f_{1}(0,1)],\ldots ,[f_{n}(0,1)])\text{ if and only
if}
\end{equation*}%
\begin{equation*}
\{\langle x,y\rangle \in M^{2}:\mathcal{M}\models R(f_{1}(x,y),\ldots
,f_{n}(x,y))\}\in \mathcal{U}^{2}.
\end{equation*}%
To see that this definition is consistent, let $f_{1}(x,y),\ldots
,f_{n}(x,y),g_{1}(x,y),\ldots ,g_{n}(x,y)$ be functions such that for all $%
1\leq i\leq n$, $[f_{i}(0,1)]=[g_{i}(0,1)]$. Now, for all $1\leq i\leq n$,
\begin{equation*}
Z_{i}=\{\langle x,y\rangle \in M^{2}:\mathcal{M}\models
f_{i}(x,y)=g_{i}(x,y)\}\in \mathcal{U}^{2}.
\end{equation*}%
Therefore $Z=Z_{1}\cap \cdots \cap Z_{n}\in \mathcal{U}^{2}$ and
\begin{equation*}
\{\langle x,y\rangle \in M^{2}:\mathcal{M}\models R(f_{1}(x,y),\ldots
,f_{n}(x,y))\}\cap Z
\end{equation*}%
\begin{equation*}
=\{\langle x,y\rangle \in M^{2}:\mathcal{M}\models R(g_{1}(x,y),\ldots
,g_{n}(x,y))\}\cap Z.
\end{equation*}%
This shows that replacing the $f_{i}$s by the $g_{i}$s does not change the
truth value of $R^{\mathcal{M}^{\ast }}([f_{1}(0,1)],\ldots ,[f_{n}(0,1)])$.

\begin{Theorems1}
(\L o\'{s} Theorem) Let $\phi (x_{1},\ldots ,x_{n})$ be an $\mathcal{L}(%
\mathcal{M})$-formula. For all $[f_{1}(0,1)],\ldots ,[f_{n}(0,1)]\in M^{\ast
}$,
\begin{equation*}
\mathcal{M}^{\ast }\models \phi ([f_{1}(0,1)],\ldots ,[f_{n}(0,1)])\text{ if
and only if}
\end{equation*}%
\begin{equation*}
\{\langle x,y\rangle \in M^{2}:\mathcal{M}\models \phi (f_{1}(x,y),\ldots
,f_{n}(x,y))\}\in \mathcal{U}^{2}.
\end{equation*}
\end{Theorems1}

\begin{proof}
We prove this theorem by structural induction on $\phi $. Without loss of
generality we may assume that $\phi $ only contains the logical connectives $%
\lnot $ and $\wedge $, and the quantifier $\exists $. It follows from the
definition of $\mathcal{M}^{\ast }$ in Subsection \ref{Warmup on two dim.
ultrapowers} and above that the theorem holds for all atomic formulae.
Suppose that the theorem holds for $\psi (x_{1},\ldots ,x_{n})$, and $\phi
(x_{1},\ldots ,x_{n})=\lnot \psi (x_{1},\ldots ,x_{n})$. Let $%
[f_{1}(0,1)],\ldots ,[f_{n}(0,1)]\in M^{\ast }$. Now,
\begin{equation*}
\mathcal{M}^{\ast }\models \phi ([f_{1}(0,1)],\ldots ,[f_{n}(0,1)])\text{ if
and only if}
\end{equation*}%
\begin{equation*}
\lnot \left( \mathcal{M}^{\ast }\models \psi ([f_{1}(0,1)],\ldots
,[f_{n}(0,1)])\right) \text{ if and only if}
\end{equation*}%
\begin{equation*}
\{\langle x,y\rangle \in M^{2}:\mathcal{M}\models \psi (f_{1}(x,y),\ldots
,f_{n}(x,y))\}\notin \mathcal{U}^{2}\text{ if and only if}
\end{equation*}%
\begin{equation*}
M^{2}\backslash \{\langle x,y\rangle \in M^{2}:\mathcal{M}\models \psi
(f_{1}(x,y),\ldots ,f_{n}(x,y))\}\in \mathcal{U}^{2}\text{ if and only if}
\end{equation*}%
\begin{equation*}
\{\langle x,y\rangle \in M^{2}:\mathcal{M}\models \lnot \psi
(f_{1}(x,y),\ldots ,f_{n}(x,y))\}\in \mathcal{U}^{2}.
\end{equation*}%
Suppose that the theorem holds for $\psi (x_{1},\ldots ,x_{n})$ and $\theta
(x_{1},\ldots ,x_{n})$, and $\phi (x_{1},\ldots ,x_{n})=\psi (x_{1},\ldots
,x_{n})\wedge \theta (x_{1},\ldots ,x_{n})$. Let $[f_{1}(0,1)],\ldots
,[f_{n}(0,1)]\in M^{\ast }$. Let
\begin{equation*}
Z_{1}=\{\langle x,y\rangle \in M^{2}:\mathcal{M}\models \psi
(f_{1}(x,y),\ldots ,f_{n}(x,y))\},
\end{equation*}%
\begin{equation*}
Z_{2}=\{\langle x,y\rangle \in M^{2}:\mathcal{M}\models \theta
(f_{1}(x,y),\ldots ,f_{n}(x,y))\},
\end{equation*}%
\begin{equation*}
Z=\{\langle x,y\rangle \in M^{2}:\mathcal{M}\models \phi (f_{1}(x,y),\ldots
,f_{n}(x,y))\}.
\end{equation*}%
Note that $Z=Z_{1}\cap Z_{2}$. Now,
\begin{equation*}
\mathcal{M}^{\ast }\models \phi ([f_{1}(0,1)],\ldots ,[f_{n}(0,1)])\text{ if
and only if}
\end{equation*}%
\begin{equation*}
\mathcal{M}^{\ast }\models \psi ([f_{1}(0,1)],\ldots ,[f_{n}(0,1)])\text{
and }\mathcal{M}^{\ast }\models \theta ([f_{1}(0,1)],\ldots ,[f_{n}(0,1)])%
\text{ if and only if}
\end{equation*}%
\begin{equation*}
Z_{1}\in \mathcal{U}^{2}\text{ and }Z_{2}\in \mathcal{U}^{2}\text{ if and
only if }Z\in \mathcal{U}^{2}.
\end{equation*}%
Suppose that the theorem holds for $\psi (w,x_{1},\ldots ,x_{n})$, and $\phi
(x_{1},\ldots ,x_{n})=\exists w\psi (w,x_{1},\ldots ,x_{n})$. Let $%
[f_{1}(0,1)],\ldots ,[f_{n}(0,1)]\in M^{\ast }$. Now, if $\mathcal{M}^{\ast
}\models \phi ([f_{1}(0,1)],\ldots ,[f_{n}(0,1)])$, then there exists $%
[h(0,1)]\in M^{\ast }$ such that
\begin{equation*}
\mathcal{M}^{\ast }\models \psi ([h(0,1)],[f_{1}(0,1)],\ldots ,[f_{n}(0,1)]).
\end{equation*}%
Therefore
\begin{equation*}
\{\langle x,y\rangle \in M^{2}:\mathcal{M}\models \psi
(h(x,y),f_{1}(x,y),\ldots ,f_{n}(x,y))\}\in \mathcal{U}^{2}.
\end{equation*}%
And so
\begin{equation*}
\{\langle x,y\rangle \in M^{2}:\mathcal{M}\models \exists w\psi
(w,f_{1}(x,y),\ldots ,f_{n}(x,y))\}\in \mathcal{U}^{2}.
\end{equation*}%
We are left to show the converse. Since $\mathcal{M}$ has definable Skolem
functions, there is an $\mathcal{M}$-definable function $g(x_{1},\ldots
,x_{n})$ such that
\begin{equation*}
\mathcal{M}\models \forall x_{1}\cdots \forall x_{n}(\exists w\psi
(w,x_{1},\ldots ,x_{n})\rightarrow \psi (g(x_{1},\ldots ,x_{n}),x_{1},\ldots
,x_{n})).
\end{equation*}%
Therefore, if
\begin{equation*}
\{\langle x,y\rangle \in M^{2}:\mathcal{M}\models \phi (f_{1}(x,y),\ldots
,f_{n}(x,y))\}\in \mathcal{U}^{2},
\end{equation*}%
\begin{equation*}
\text{then }\{\langle x,y\rangle \in M^{2}:\mathcal{M}\models \psi
(g(f_{1}(x,y),\ldots ,f_{n}(x,y)),f_{1}(x,y),\ldots ,f_{n}(x,y))\}\in
\mathcal{U}^{2}.
\end{equation*}%
On the other hand, since $h(x,y)=g(f_{1}(x,y),\ldots ,f_{n}(x,y))$ is $%
\mathcal{M}$-definable, we have
\begin{equation*}
\mathcal{M}^{\ast }\models \psi ([h(0,1)],[f_{1}(0,1)],\ldots ,[f_{n}(0,1)]),
\end{equation*}%
\begin{equation*}
\text{and so }\mathcal{M}^{\ast }\models \phi ([f_{1}(0,1)],\ldots
,[f_{n}(0,1)]).
\end{equation*}%
The theorem now follows by induction.\hfill $\square $\medskip
\end{proof}

\noindent \textbf{(b).} We need to show that $\{[0],[1]\}$ forms a set of
order indiscernibles in $\mathcal{M}^{\ast }$ over $\mathcal{M}$. Let $\phi
(x_{0},\ldots ,x_{n})$ be an $\mathcal{L}(\mathcal{M})$-formula. Let $%
m_{1},\ldots ,m_{n}\in M$. For each $1\leq i\leq n$, $m_{i}$ is represented
in $\mathcal{M}^{\ast }$ by the constant function $h_{i}(x,y)=m_{i}$. Let
\begin{equation*}
Z=\{w\in M:\mathcal{M}\models \phi (w,m_{1},\ldots ,m_{n})\}.
\end{equation*}%
Note that
\begin{equation*}
\{m\in M:(M\times Z)|m\in \mathcal{U}\}=\left\{
\begin{array}{ll}
\emptyset & \text{if }Z\notin \mathcal{U} \\
M & \text{if }Z\in \mathcal{U}%
\end{array}%
\right.
\end{equation*}%
\begin{equation*}
\text{and }\{m\in M:(Z\times M)|m\in \mathcal{U}\}=Z.
\end{equation*}%
Therefore, $M\times Z\in \mathcal{U}^{2}$ if and only if $Z\times M\in
\mathcal{U}^{2}$. Now, using part (a),
\begin{equation*}
\mathcal{M}^{\ast }\models \phi ([0],[h_{1}(0,1)],\ldots ,[h_{n}(0,1)])\text{
if and only if}
\end{equation*}%
\begin{equation*}
\{\langle x,y\rangle \in M^{2}:\mathcal{M}\models \phi (x,m_{1},\ldots
,m_{n})\}=Z\times M\in \mathcal{U}^{2}\text{ if and only if}
\end{equation*}%
\begin{equation*}
\{\langle x,y\rangle \in M^{2}:\mathcal{M}\models \phi (y,m_{1},\ldots
,m_{n})\}=M\times Z\in \mathcal{U}^{2}\text{ if and only if}
\end{equation*}%
\begin{equation*}
\mathcal{M}^{\ast }\models \phi ([1],[h_{1}(0,1)],\ldots ,[h_{n}(0,1)]).
\end{equation*}%
This shows that $\{[0],[1]\}$ forms a set of order indiscernibles in $%
\mathcal{M}^{\ast }$ over $\mathcal{M}$.

\noindent \textbf{3.9. Exercise.} Let $\mathcal{U}$ be an $\mathcal{M}$%
-amenable ultrafilter. Let $(I,<)$ be an ordered set disjoint from $M$. Let $%
\mathcal{M}^{\ast }=\mathrm{Ult}(\mathcal{M},\mathcal{U},I)$ (defined on
p.16-17). We need to define the interpretations of the relation symbols in $%
\mathcal{L}(\mathcal{M})$ in $\mathcal{M}^{\ast }$. Let $R(x_{1},\ldots
,x_{n})$ be an $n$-ary relation symbol in $\mathcal{L}(\mathcal{M})$. For
all $[f_{1}(I_{1})],\ldots ,[f_{n}(I_{n})]\in M^{\ast }$, define
\begin{equation*}
R^{\mathcal{M}^{\ast }}([f_{1}(I_{1})],\ldots ,[f_{n}(I_{n})])\text{ if and
only if}
\end{equation*}%
\begin{equation*}
\{u\in M^{I_{0}}:\mathcal{M}\models R(f_{1}(I_{1})[u],\ldots
,f_{n}(I_{n})[u])\}\in \mathcal{U}^{I_{0}}\text{ where }I_{0}=\bigcup_{1\leq
j\leq n}I_{j}.
\end{equation*}%
Let $\varepsilon :M\longrightarrow M^{\ast }$ be the map defined by $%
m\mapsto \lbrack f_{m}(I_{0})]$ where $f_{m}$ is the constant map with value
$m$ and $I_{0}\subseteq I$ is finite. We need to verify that $\varepsilon $
is an isomorphism between $\mathcal{M}$ and $\varepsilon (\mathcal{M})$. Let
$m_{1},m_{2}\in M$. We have
\begin{equation*}
\varepsilon (m_{1})=\varepsilon (m_{2})\text{ if and only if}
\end{equation*}%
\begin{equation*}
\{u\in M^{I_{2}}:\mathcal{M}\models f_{m_{1}}(I_{0})=f_{m_{2}}(I_{1})\}\in
\mathcal{U}^{I_{2}}
\end{equation*}%
\begin{equation*}
\text{ where }I_{0},I_{1}\subseteq I\text{ are finite and }I_{2}=I_{0}\cup
I_{1}
\end{equation*}%
\begin{equation*}
\text{if and only if }m_{1}=m_{2}.
\end{equation*}%
This shows that $\varepsilon $ is a well-defined injection. Let $%
R(x_{1},\ldots ,x_{n})$ be an $n$-ary relation symbol in $\mathcal{L}(%
\mathcal{M})$. Let $m_{1},\ldots ,m_{n}\in M$. We have
\begin{equation*}
\mathcal{M}^{\ast }\models R(\varepsilon (m_{1}),\ldots ,\varepsilon (m_{n}))%
\text{ if and only if}
\end{equation*}%
\begin{equation*}
\{u\in M^{I_{0}}:\mathcal{M}\models R(f_{m_{1}}(I_{1})[u],\ldots
,f_{m_{n}}(I_{n})[u])\}\in \mathcal{U}^{I_{0}}
\end{equation*}%
\begin{equation*}
\text{where }I_{1},\ldots ,I_{n}\subseteq I\text{ are finite and }%
I_{0}=\bigcup_{1\leq j\leq n}I_{j}
\end{equation*}%
\begin{equation*}
\text{if and only if }\mathcal{M}\models R(m_{1},\ldots ,m_{n}).
\end{equation*}%
Let $g(x_{1},\ldots ,x_{n})$ be an $n$-ary function symbol in $\mathcal{L}(%
\mathcal{M})$. Let $m_{1},\ldots ,m_{n}\in M$. Let $m_{0}\in M$ be such that
\begin{equation*}
\mathcal{M}\models m_{0}=g(m_{1},\ldots ,m_{n}).
\end{equation*}%
Now,
\begin{equation*}
g^{\mathcal{M}^{\ast }}(\varepsilon (m_{1}),\ldots ,\varepsilon (m_{n}))=g^{%
\mathcal{M}^{\ast }}([f_{m_{1}}(I_{1})],\ldots
,[f_{m_{n}}(I_{n})])=[h(I_{0})],
\end{equation*}%
where $I_{1},\ldots ,I_{n}\subseteq I$ are finite and $I_{0}=\bigcup_{1\leq
j\leq n}I_{j}$, and $h$ is a parametrically definable function such that for
all $u\in M^{I_{0}}$,
\begin{equation*}
h(I_{0})[u]=g^{\mathcal{M}}(f_{m_{1}}(I_{1}),\ldots ,f_{m_{n}}(I_{n}))=g^{%
\mathcal{M}}(m_{1},\ldots ,m_{n})=m_{0}.
\end{equation*}%
Therefore $h(I_{0})=f_{m_{0}}(I_{0})=\varepsilon (m_{0})$. This shows that $%
\varepsilon $ is an isomorphism.

\noindent \textbf{3.10. Exercise.} Let $I_{0}\subseteq I$ be finite. We show
that if $I_{0}=I_{1}\cup I_{2}$ with $\max I_{1}<\min I_{2}$, then for all $%
X\subseteq M^{I_{0}}$,
\begin{equation*}
X\in \mathcal{U}^{I_{0}}\text{ if and only if }\{s\in M^{I_{1}}:X|s\in
\mathcal{U}^{I_{2}}\}\in \mathcal{U}^{I_{1}},
\end{equation*}%
by induction on the size of $I_{1}$. When $|I_{1}|=1$ the result follows
immediately from the definition of $\mathcal{U}^{I_{0}}$. Suppose that the
result holds for all finite $J_{0}=J_{1}\cup J_{2}\subseteq I$ with $\max
J_{1}<\min J_{2}$ and $|J_{1}|=n$. Let $K_{0}\subseteq I$ and suppose $%
K_{0}=K_{1}\cup K_{2}$ with $\max K_{1}<\min K_{2}$ and $|K_{1}|=n+1$.
Suppose $K_{1}=\{i_{0}<\cdots <i_{n}\}$, and let $K_{1}^{\prime
}=\{i_{1}<\cdots <i_{n}\}$. Let $X\subseteq M^{K_{0}}$. From the definition
of $X|s$, we get
\begin{equation*}
\{s\in M^{K_{1}}:X|s\in \mathcal{U}^{K_{2}}\}\in \mathcal{U}^{K_{1}}\text{
if and only if}
\end{equation*}%
\begin{equation*}
\{q\in M^{K_{1}}:q=s\cup t\text{ where }s\in M^{\{i_{0}\}}\text{ and }t\in
M^{K_{1}^{\prime }}\text{ and }(X|s)|t\in \mathcal{U}^{K_{2}}\}\in \mathcal{U%
}^{K_{1}}.
\end{equation*}%
Therefore the definition of $\mathcal{U}^{K_{1}}$ yields
\begin{equation*}
\{s\in M^{K_{1}}:X|s\in \mathcal{U}^{K_{2}}\}\in \mathcal{U}^{K_{1}}\text{
if and only if}
\end{equation*}%
\begin{equation*}
\{s\in M^{\{i_{0}\}}:\{t\in M^{K_{1}^{\prime }}:(X|s)|t\in \mathcal{U}%
^{K_{2}}\}\in \mathcal{U}^{K_{1}^{\prime }}\}\in \mathcal{U}^{\{i_{0}\}}.
\end{equation*}%
Which, by the induction hypothesis, gives
\begin{equation*}
\{s\in M^{K_{1}}:X|s\in \mathcal{U}^{K_{2}}\}\in \mathcal{U}^{K_{1}}
\end{equation*}%
\begin{equation*}
\text{if and only if }\{s\in M^{\{i_{0}\}}:X|s\in \mathcal{U}^{K_{1}^{\prime
}\cup K_{2}}\}\in \mathcal{U}^{\{i_{0}\}}.
\end{equation*}%
And finally, the definition of $\mathcal{U}^{K_{0}}$ gives
\begin{equation*}
\{s\in M^{K_{1}}:X|s\in \mathcal{U}^{K_{2}}\}\in \mathcal{U}^{K_{1}}\text{
if and only if }X\in \mathcal{U}^{K_{0}}.
\end{equation*}%
This completes the proof.

\noindent \textbf{3.12. Exercise.} Let $I_{2}\subseteq I$ be finite with $%
I_{2}=I_{0}\cup I_{1}$. Let $J\subseteq I$ be finite with $I_{2}\subseteq J$%
. Let $f(I_{0})$ and $g(I_{1})$ be generalized terms. Let
\begin{equation*}
Z_{1}=\{u\in M^{I_{2}}:f(I_{0})[u]=g(I_{1})[u]\}\text{ and }Z_{2}=\{u\in
M^{J}:f(I_{0})[u]=g(I_{1})[u]\}.
\end{equation*}%
We need to show that $Z_{1}\in \mathcal{U}^{I_{2}}$ if and only if $Z_{2}\in
\mathcal{U}^{J}$. But this follows immediately from Lemma 3.11, since
\begin{equation*}
Z_{2}=\{s\cup t\in M^{J}:s\in Z_{1}\wedge t\in M^{J\backslash I_{2}}\}.
\end{equation*}

\noindent \textbf{3.13. Exercise.}

\begin{Theorems1}
(\L o\'{s} Theorem) Let $\phi (x_{1},\ldots ,x_{n})$ be an $\mathcal{L}(%
\mathcal{M})$-formula. For all $[f_{1}(I_{1})],\ldots ,[f_{n}(I_{n})]\in
M^{\ast }$,
\begin{equation*}
\mathcal{M}^{\ast }\models \phi ([f_{1}(I_{1})],\ldots ,[f_{n}(I_{n})])\text{
if and only if}
\end{equation*}%
\begin{equation*}
\{u\in M^{I_{0}}:\mathcal{M}\models \phi (f_{1}(I_{1})[u],\ldots
,f_{n}(I_{n})[u])\}\in \mathcal{U}^{I_{0}}\text{ where }I_{0}=\bigcup_{1\leq
j\leq n}I_{j}.
\end{equation*}
\end{Theorems1}

\begin{proof}
We prove this theorem by structural induction on $\phi $. Without loss of
generality we may assume that $\phi $ only contains the logical connectives $%
\lnot $ and $\wedge $, and the quantifier $\exists $. It follows from the
definition of $\mathcal{M}^{\ast }$ in Subsection \ref{Dimensional
Ultrapowers} and Exercise 3.9 that the theorem holds for all atomic
formulae. As was the case with the \L o\'{s} Theorem proved in Exercise 3.8,
the inductive steps where $\phi =\lnot \psi $ and $\phi =\psi \wedge \theta $
follow from the basic properties of the ultrafilter $\mathcal{U}^{I_{0}}$.
Again, the nontrivial inductive step involves dealing with the quantifier $%
\exists $. Suppose that the theorem holds for $\psi (w,x_{1},\ldots ,x_{n})$%
, and $\phi (x_{1},\ldots ,x_{n})=\exists w\psi (w,x_{1},\ldots ,x_{n})$.
Let $[f_{1}(I_{1})],\ldots ,[f_{n}(I_{n})]\in M^{\ast }$. If
\begin{equation*}
\mathcal{M}^{\ast }\models \phi ([f_{1}(I_{1})],\ldots ,[f_{n}(I_{n})]),
\end{equation*}%
then there exists $[g(J)]\in M^{\ast }$ such that
\begin{equation*}
\mathcal{M}^{\ast }\models \psi ([g(J)],[f_{1}(I_{1})],\ldots
,[f_{n}(I_{n})]).
\end{equation*}%
Therefore, by the induction hypothesis,
\begin{equation*}
Z_{1}=\{u\in M^{J^{\prime }}:\mathcal{M}\models \psi
(g(J)[u],f_{1}(I_{1})[u],\ldots ,f_{n}(I_{n})[u])\}\in \mathcal{U}%
^{J^{\prime }}
\end{equation*}%
\begin{equation*}
\text{where }J^{\prime }=J\cup \bigcup_{1\leq j\leq n}I_{j}.
\end{equation*}%
Let $I_{0}=\bigcup_{1\leq j\leq n}I_{j}$. Let
\begin{equation*}
Z_{2}=\{u\in M^{I_{0}}:\mathcal{M}\models \exists w\psi
(w,f_{1}(I_{1})[u],\ldots ,f_{n}(I_{n})[u])\}\text{ and}
\end{equation*}%
\begin{equation*}
Z_{3}=\{u\in M^{J^{\prime }}:\mathcal{M}\models \exists w\psi
(w,f_{1}(I_{1})[u],\ldots ,f_{n}(I_{n})[u])\}.
\end{equation*}%
Now, $Z_{1}\subseteq Z_{3}$, and so $Z_{3}\in \mathcal{U}^{J^{\prime }}$.
And
\begin{equation*}
Z_{3}=\{s\cup t\in M^{J^{\prime }}:s\in Z_{2}\wedge t\in M^{J^{\prime
}\backslash I_{0}}\},
\end{equation*}%
so, by Lemma 3.11, $Z_{2}\in \mathcal{U}^{I_{0}}$.\newline
We are left to show the converse. Since $\mathcal{M}$ has definable Skolem
functions, there exists a function $g(x_{1},\ldots ,x_{n})$ such that
\begin{equation*}
\mathcal{M}\models \forall x_{1}\cdots \forall x_{n}(\exists w\psi
(w,x_{1},\ldots ,x_{n})\rightarrow \psi (g(x_{1},\ldots ,x_{n}),x_{1},\ldots
,x_{n})).
\end{equation*}%
Let $I_{0}=\bigcup_{1\leq j\leq n}I_{j}$. Now, if
\begin{equation*}
\{u\in M^{I_{0}}:\mathcal{M}\models \phi (f_{1}(I_{1})[u],\ldots
,f_{n}(I_{n})[u])\}\in \mathcal{U}^{I_{0}},
\end{equation*}%
\begin{equation*}
\text{then }\{u\in M^{I_{0}}:\mathcal{M}\models \psi
(g(f_{1}(I_{1})[u],\ldots ,f_{n}(I_{n})[u]),f_{1}(I_{1})[u],\ldots
,f_{n}(I_{n})[u])\}\in \mathcal{U}^{I_{0}}.
\end{equation*}%
Let $h(I_{0})=g(f_{1}(I_{1}),\ldots ,f_{n}(I_{n}))$. Therefore
\begin{equation*}
\{u\in M^{I_{0}}:\mathcal{M}\models \psi (h(I_{0})[u],f_{1}(I_{1})[u],\ldots
,f_{n}(I_{n})[u])\}\in \mathcal{U}^{I_{0}}.
\end{equation*}%
And so, by the induction hypothesis,
\begin{equation*}
\mathcal{M}^{\ast }\models \psi ([h(I_{0})],[f_{1}(I_{1})],\ldots
,[f_{n}(I_{n})]),
\end{equation*}%
\begin{equation*}
\text{which means }\mathcal{M}^{\ast }\models \phi ([f_{1}(I_{1})],\ldots
,[f_{n}(I_{n})]).
\end{equation*}%
The theorem now follows by induction.\hfill $\square $\medskip
\end{proof}

We now turn to showing that the embedding $\varepsilon :\mathcal{M}%
\longrightarrow \mathcal{M}^{\ast }$ defined in Subsection \ref{Dimensional
Ultrapowers} is elementary. As we did in the solution to Exercise 3.9, if $%
I_{0}\subseteq I$ is finite and $m\in M$, then we write $f_{m}(I_{0})$ for
the constant function with value $m$. Let $\phi (x_{1},\ldots ,x_{n})$ be an
$\mathcal{L}(\mathcal{M})$-formula. Let $I_{0}\subseteq I$ be finite. Note
that, for all $m_{1},\ldots ,m_{n}\in M$,
\begin{equation*}
\{u\in M^{I_{0}}:\mathcal{M}\models \phi (f_{m_{1}}(I_{0})[u],\ldots
,f_{m_{n}}(I_{0})[u])\}=\left\{
\begin{array}{ll}
M^{I_{0}} & \text{if }\mathcal{M}\models \phi (m_{1},\ldots ,m_{n}) \\
\emptyset & \text{otherwise}%
\end{array}%
\right.
\end{equation*}%
Therefore the \L o\'{s} Theorem proved above yields:
\begin{equation*}
\mathcal{M}^{\ast }\models \phi (\varepsilon (m_{1}),\ldots ,\varepsilon
(m_{n}))\text{ if and only if}
\end{equation*}%
\begin{equation*}
\{u\in M^{I_{0}}:\mathcal{M}\models \phi (f_{m_{1}}(I_{0})[u],\ldots
,f_{m_{n}}(I_{0})[u])\}\in \mathcal{U}^{I_{0}}\text{ if and only if}
\end{equation*}%
\begin{equation*}
\mathcal{M}\models \phi (m_{1},\ldots ,m_{n}).
\end{equation*}

\noindent \textbf{3.14. Exercise.} We need to show that $(I,<)$ is a set of
order indiscernibles over $\mathcal{M}$. Let $\phi (y_{1},\ldots
,y_{n},x_{0},\ldots ,x_{k-1})$ be an $\mathcal{L}(\mathcal{M})$-formula. Let
$m_{1},\ldots ,m_{n}\in M$, and let $i_{0}<\cdots <i_{k-1},j_{0}<\cdots
<j_{k-1}$ be in $I$. Let $I_{0}=\{i_{0}<\cdots <i_{k-1}\}$ and let $%
J_{0}=\{j_{0}<\cdots <j_{k-1}\}$. Now,
\begin{equation*}
\mathcal{M}^{\ast }\models \phi (\varepsilon (m_{1}),\ldots ,\varepsilon
(m_{n}),[i_{0}],\ldots ,[i_{k-1}])\text{ if and only if}
\end{equation*}%
\begin{equation*}
\{u\in M^{I_{0}}:\mathcal{M}\models \phi (m_{1},\ldots
,m_{n},i_{0}[u],\ldots ,i_{k-1}[u])\}\in \mathcal{U}^{I_{0}}\text{ if and
only if}
\end{equation*}%
\begin{equation*}
\{\langle x_{0},\ldots ,x_{k-1}\rangle \in M^{k}:\mathcal{M}\models \phi
(m_{1},\ldots ,m_{n},x_{0},\ldots ,x_{k-1})\}\in \mathcal{U}^{k}\text{ if
and only if}
\end{equation*}%
\begin{equation*}
\{u\in M^{J_{0}}:\mathcal{M}\models \phi (m_{1},\ldots
,m_{n},j_{0}[u],\ldots ,j_{k-1}[u])\}\in \mathcal{U}^{J_{0}}\text{ if and
only if}
\end{equation*}%
\begin{equation*}
\mathcal{M}^{\ast }\models \phi (\varepsilon (m_{1}),\ldots ,\varepsilon
(m_{n}),[j_{0}],\ldots ,[j_{k-1}]),
\end{equation*}%
where the first and last equivalences follow from the \L o\'{s} Theorem, and
the middle two equivalences follow from the definition of the dimensional
ultrapower.

\end{document}